\documentclass[12pt]{article}
\usepackage{fancyhdr}
\usepackage{amssymb}
\usepackage{amsmath}
\usepackage{amsthm}
\usepackage{graphicx}
\usepackage{psfrag}
\usepackage{graphics}
\usepackage{amsthm}
\usepackage{hyperref,enumitem}
\usepackage[dvipsnames]{xcolor}
\usepackage[english]{babel}

\newcommand{\UD}{\mathcal{UD}_{P}}
\newcommand{\LTC}{\mathcal{LTC}}

\newcommand{\cI}{\mathcal{I}}

\newcommand{\cR}{\mathcal{R}}
\newcommand{\cV}{\mathcal{V}}

\newcommand{\N}{\mathbb N}

\newcommand{\Homeo}{\mathrm{Homeo}}
\newcommand{\Int}{\mathrm{int}}
\newcommand{\var}{\mathrm{var}}

\newcommand{\R}{\mathbb{R}}

\newcommand{\Z}{\mathbb Z}

\newcommand{\Diff}{\mathrm{Diff}}

\newcommand{\eps}{\varepsilon}

\def\Fix{{\mathrm {Fix}}}

\theoremstyle{theorem}

\newtheorem{thm}{Theorem}[section]
\newtheorem{thmintro}{Theorem}

\newtheorem{prop}[thm]{Proposition}

\newtheorem{cor}[thm]{Corollary}

\newtheorem{lem}[thm]{Lemma}
\newtheorem{claim}{Claim}

\theoremstyle{definition}
\newtheorem{defn}[thm]{Definition}
\theoremstyle{remark}
\newtheorem{rem}[thm]{Remark}

\newcommand{\id}{\mathrm{id}}

\begin{document}

\date{}
\vspace{-1cm}

\date{}
\author{H\'el\`ene Eynard-Bontemps \,\, \& \,\, Andr\'es Navas}

\title {All, most, some?\\
On diffeomorphisms of the interval that are distorted \\
and/or conjugate to powers of themselves}

\maketitle

\noindent{\bf Abstract:} We study the problem of conjugating a diffeomorphism of the interval to (positive) powers of 
itself.  Although this is always possible for homeomorphisms, the smooth setting is rather interesting. Besides the 
obvious obstruction given by hyperbolic fixed points, several other aspects need to be considered. As concrete 
results we show that, in class $C^1$, if we restrict to the (closed) subset of diffeomorphisms having only parabolic 
fixed points, then the set of diffeomorphisms that are conjugate to their powers is dense, but its complement is 
generic. In higher regularity, however, the complementary set contains an open and dense set.  The 
text is complemented with several remarks and results concerning distortion 
elements of the group of diffeomorphisms of the interval in several regularities.

\vspace{0.2cm}

\noindent{\bf Keywords:} Diffeomorphism, conjugacy, Mather invariant, asymptotic variation, distortion element.

\vspace{0.2cm}

\noindent{\bf MCS 2020:} 
37C10, 
37C15, 
37C85, 
37E05, 
57M60. 

\section{Introduction}

This work deals with a very concrete problem: 

\vspace{0.4cm}

\noindent{\bf General Question:} What are the diffeomorphisms $f$ of the interval that are conjugated to some (all) of their 
positive powers $f^k$~? \medskip

\vspace{0.4cm}

Obviously, we are thinking about ``nontrivial'' powers, meaning that $k \geq 2$. 

There are many motivations for this question. For instance, it naturally arises when studying actions of the affine group. 
Another important source of inspiration comes from the notion of {\em distortion element} introduced by Gromov, which 
has been extensively exploited in dynamical contexts over the last 20 years (see for instance 
\cite{CF06,FH06,Av08,Mi14,LM18,Mi18,MR18,GL19,HM19,Na21,Ro21,BS22,DE22}).

Recall that an element $f$ of a finitely-generated group $G$ is said to be {\em distorted} if for some (equivalently, every) finite system of generators, 
the corresponding word-length of the powers $f^n$ grows sublinearly in $n$. More generally, an element $f$ of a general group $G$ 
is a {\em distortion element} if it is distorted inside some finitely generated subgroup. It turns out that if $f$ is conjugated to a 
nontrivial power $f^k$ in $G$, then it is a distortion element. Indeed, from $f^k = h f h^{-1}$ one easily concludes that 
$$f^{k^m} = h^m f h^{-m},$$
which shows that the word-length of $f^{k^m}$ is at most of the order $m$, and easily allows showing that the world-length 
of $f^n$ is at most of the order $\log_k (n)$ (see \cite{Di23} for further discussion on this argument). Note, however, that this 
is far from being the only procedure to produce distortion elements inside groups;  see for instance 
\cite[Appendix]{Gr81} for relevant examples on this.

In this work we only deal with orientation-preserving maps (compare \cite{GOS09}). 
With this restriction, an elementary argument shows that every interval homeomorphism is conjugated to each of its nontrivial 
powers (cf. \cite[p. 348]{Gh01}). In particular, every such homeomorphism is a distortion element of $\Homeo_+([0,1])$. In 
regularity $C^1$ and higher, there is an obvious obstruction to distortion, namely all fixed points 
must be parabolic, otherwise the logarithm of the derivative of the powers would grow linearly (cf. \cite{PS04}). We will thus restrict 
to the subspace of diffeomorphisms having only parabolic fixed points, that we denote $\Diff_P^{\ell} ([0,1])$ (here, ``$P\,$'' stands 
for ``parabolic'', and $\ell$ refers to the regularity class in consideration). It is worth noting that this is a closed subspace  
of the corresponding diffeomorphism group. Indeed, its complementary set is open, as it directly follows from the 
persistence of hyperbolic fixed points under $C^1$-small perturbations. In particular, for each integer $\ell \geq 1$, 
the group $\Diff_P^{\ell}([0,1])$ is a complete metric space when endowed with the natural metric (recalled below). 
As a consequence, it has the Baire property, and thus one is allowed to look for generic properties inside. Also, 
though not being a metric space, $\Diff^{\infty}_{P}([0,1])$ is a Baire space as well when endowed with its natural topology.

For simplicity, if $f$ is a distortion element of 
$\Diff_+^{\ell}([0,1])$, then sometimes it will be just said to be {\em $C^{\ell}$-distorted}.
Amazingly, the answer to the following question is unknown:

\vspace{0.4cm}
\noindent{\bf Question 1.} 
Is every element in $\Diff^1_P ([0,1])$ a distortion element of $\Diff_+^1([0,1])$?

\vspace{0.4cm}

According to \cite{Co20}, an affirmative answer to this question would imply that another notion of distortion introduced 
by Rosendal in \cite{Ro21} would coincide with the one considered here in the particular case of the group $\Diff^1_+([0,1])$. 

Although we do not solve this problem, we show two results concerning it, the first one in the positive direction 
and the second one somewhat in the negative.

\begin{thmintro}
\label{t:main}
Inside the space $\Diff^1_P([0,1])$, the $C^{\infty}$ diffeomorphisms that are $C^\infty$-conjugated to each of their positive powers 
form a dense subset.
\end{thmintro}

\begin{thmintro}
\label{t:generic}
Inside the space $\Diff^1_P([0,1])$, the property of \emph{not} being conjugated to any of its nontrivial powers is generic. 
\end{thmintro}

The answer to Question 1 might nevertheless be positive, but if so, most of the elements of $\Diff^1_P([0,1])$ would 
be $C^1$-distorted for much less obvious reasons than just being conjugate to a nontrivial power. An interesting family 
of such examples is given in \cite{Na21}; see also \S \ref{s:remarks-fin} for another family of a very different nature.

 \medskip

The situation is rather different in higher regularity, as it is shown by the next result. 

\begin{thmintro} 
\label{t:higher}
For each $\ell \geq 2$ (integer or infinite), the subset of $\mathrm{Diff}^{\ell}_P ([0,1])$ made of elements that are neither 
$C^{\ell}$-distorted nor $C^1$-conjugated to some nontrivial power of themselves contains an open and dense set, yet it is not open.
\end{thmintro}

Roughly, the open and dense set alluded to in this theorem is a large subset where a natural obstruction to $C^2$-distortion arises, 
namely the nonvanishing of the so-called {\em asymptotic variation}.\footnote{
This was first introduced in \cite{Na23} under the somewhat misleading name of \emph{asymptotic distortion} 
and referred to in this way in \cite{EN21}, where it is further developed.} 
For a $C^2$ diffeomorphism $f$ (and more generally, for a $C^{1+bv}$ diffeomorphism; see below for a review on this 
regularity class), this is defined as the limit 
$$\lim_{n \to \infty} \frac{\var (\log Df^n)}{n},$$
and denoted $V_{\infty} (f)$. This limit exists by the subadditivity of $f  \mapsto \var (\log Df)$. This subadditivity also implies that, 
if $f$ is $C^{1+bv}$-distorted, then $V_{\infty} (f) = 0$. (See \S \ref{s:higher} for more details.)
The following question is then the analog to Question~1 in higher regularity.  
For $\ell \geq 2$ or $\ell = 1+ bv$, we denote by $\Diff^{\ell}_V ([0,1])$ the subset of $\Diff^{\ell}_+ ([0,1])$ made of the 
diffeomorphisms with vanishing asymptotic variation.

\medskip

\noindent{\bf Question 2.} Are all diffeomorphisms $f \in \Diff^{\ell}_V ([0,1])$ necessarily $C^{\ell}$-distorted~? 

\medskip

We give a partial positive answer to this question in class $C^2$.

\begin{thmintro} 
\label{t:higher2}
The set of diffeomorphisms that are $C^2$-conjugated to all their nontrivial powers is a dense 
subset of  $\Diff^{2}_V ([0,1])$.  
\end{thmintro}

We will see (cf. Remark \ref{r:last}) that this result does not extend to higher regularity.

\medskip

Note that, since we deal with dense/generic sets, the dichotomy in all our results is very strong: they involve 
diffeomorphisms that are conjugate to all/none of their nontrivial powers. This leaves open the following interesting question:  

\vspace{0.4cm}

\noindent{\bf Question 3.} Does there exist  $\ell \ge 1$ and $f \in \Diff^{\ell}_P ([0,1])$ that is conjugate to some nontrivial power 
but not to another one~? For example, does there exist such a diffeomorphism that is conjugate to its square but not to its cube~? 
Here we are thinking about $C^1$ or $C^{\ell}$ conjugacy (the question makes sense in both categories and it may be the case 
that it has different answers depending on the regularity). 

\vspace{0.4cm}

Though we are primarily concerned with $C^1$ diffeomorphisms, part of the above statements and/or their proofs involve diffeomorphisms of higher 
regularity. For these, the $C^1$ conjugacy classes are much better understood, mainly thanks to classical works by \cite{szekeres,kopell,Ma?,yoccoz}. 
We thus start this article with a background section (cf. \S \ref{s:back}) recalling these results. (The reader who is already familiar with 
these notions can directly go to the proofs, and the one who is not can proceed similarly and go back to the background section whenever necessary.) 
The proof of Theorem \ref{t:main} is next given in \S \ref{s:distorted}. For Theorem \ref{t:generic}, two independent proofs are given in \S \ref{s:generic-short} 
and \S \ref{s:generic}, respectively. Theorem \ref{t:higher} is proved in \S \ref{s:higher}. Finally, Theorem \ref{t:higher2} is proved in \S \ref{s:ln}.

\vspace{0.4cm}

\noindent{\bf Some notation.} 
We close this introduction by fixing some notation and vocabulary that will be used throughout this work. Given an integer 
$\ell \geq 0$ and $f,g$ in $\Diff^{\ell}_+([0,1])$, we let 
$$\|f-g\|_{\ell} := \max_{x \in [0,1]} \big| D^{(\ell)} f(x) -D^{(\ell)} g(x) \big| .$$ 
This defines a distance in $\Diff^{\ell}_+ ([0,1])$ that determines its canonical topology. 
For circle diffeomorphisms, the distance we consider is given by 
$$\|f-g\|_{\ell} := \max_{x \in \mathbb{S}^1} \big| f(x)- g (x) \big| + \max_{x \in \mathbb{S}^1} \big| D^{(\ell)} f(x) -D^{(\ell)} g(x) \big| .$$ 
Say that $f,g$ are $\varepsilon$-close in $C^{\ell}$ norm if $\| f - g \|_{\ell} \leq \varepsilon$. 
Similar notation will be used for restrictions of diffeomorphisms to subintervals.

We also consider the group $\Diff^{1+bv}_+([0,1])$ of $C^1$ diffeomorphisms 
whose derivative (equivalently, the logarithm of it) has bounded variation. 
For $f$ therein, we use the notation $V(f)$ (instead of the more standard notation $\var (\log Df)$) 
for the total variation of $ \log (Df)$: 
$$V(f) := \sup_{a_0<a_1<\ldots<a_n} \sum_{i=0}^{n-1} \big| \log Df (a_{i+1}) - \log Df (a_i) \big|.$$
Recall that $d_* (f,g) := \mathrm{var} (\log Df - \log Dg)$ defines a (complete) metric on $\Diff^{1+bv}_+([0,1])$ 
(yet it does not endow it with the structure of a topological group; cf. \cite[Appendix 4]{EN24}). 
We will also use the notation $V(f;[a,b])$ for the total variation of $\log(Df)$ either for a diffeomorphism $f$ 
of an interval $[a,b]$ or for the restriction to a subinterval $[a,b]$ of a diffeomorphism $f$ of some larger interval.

Given a diffeomorphism $f$ of $[0,1]$, we will denote by $\Fix(f)$ the set of its fixed points. 
By a \emph{component of $f$} we will refer to the closure of any connected component of $[0,1]\setminus \Fix(f)$. 
The interior of such an $I$ will be denoted by $I^{\mathrm{o}}$. 
If $f$ is not the identity, we will call a {\em fundamental domain} of $f$ any closed subinterval of the form 
$[f(x),x]$ (if $f(x)<x$) or $[x,f(x)]$ (if $f(x)>x$), where $x\in[0,1]\setminus \Fix(f)$. 
The orbit of a point $x$ under a map $f$ will be denoted $O_f (x)$. Given $l\in\N$ (which might be different from the regularity of $f$), 
we will denote by $\mathrm{Z}_f^l$ the $C^l$-centralizer of $f$, that is, the set $\{g\in\Diff^l_+([0,1]) \!:  g\circ f = f\circ g\}$.

We will denote $T_{\tau}$ the translation by $\tau$ on the real line. Finally, letters $k,n$ will denote positive integer numbers. 
The same will apply to $\ell$ yet, since it will refer to regularity classes, it can also take the value $\ell = \infty$ 
and $\ell = 1+bv$, the latter to refer to $C^1$ diffeomorphisms whose derivative has bounded variation. In 
particular, note that the statement of Theorem \ref{t:higher} includes the $C^{\infty}$ setting (where 
the topology is the one induced by the family of metrics $\|\cdot\|_{\ell}$ for finite $\ell$).


\section{Background: some classical results about interval diffeomorphisms}
\label{s:back}

Most of the results in this section can be found in \cite[Chapters IV and V]{yoccoz}. Slightly weaker extensions 
to $C^{1+bv}$ diffeomorphisms appear in \cite[Appendix 2]{EN23} (below this regularity class, the results are 
far from being true). 

Let $f$ be a $C^{\ell}$ diffeomorphism of the half-closed interval $[0,1)$ with no fixed point in the interior, with $\ell \geq 2$. 
According to Szekeres \cite{szekeres} and Kopell \cite{kopell}, there exists a unique vector field $X_f$ on $[0,1)$ such that: 

\vspace{0.1cm}

\noindent - $X_f$ is $C^1$ on $[0,1)$ (and extends continuously to $[0,1]$ by letting $X(1)=0$);

\vspace{0.1cm}

\noindent - $f$ is the time-1 map of the flow of $X_f$.

\vspace{0.1cm}

\noindent Moreover, this vector field is $C^{\ell-1}$ on $(0,1)$, and its flow coincides with the centralizer of $f$ 
in $\Diff^1_+([0,1))$.

\vspace{0.2cm}

Let now $f$ be a $C^{\ell}$ diffeomorphism of the {\em closed} interval $[0,1]$ with no fixed point in the interior, 
where $\ell \geq 2$. 
The above results apply to the restrictions of $f$ to $[0,1)$ and $(0,1]$ respectively. This provides two ``generating vector fields'' 
$X_f$ and $Y_f$ for $f$ (meaning that $f$ is the time-$1$ map of both of them), which are of class $C^1$ on $[0,1)$ and $(0,1]$, 
respectively. As in Yoccoz \cite{yoccoz}, we will denote by $f_t$ and $f^t$ the time-$t$ maps of the respective flows (so that 
$f_1=f^1=f$), which correspond to the elements of the centralizers of $f$ in $\Diff^1_+([0,1))$ and $\Diff^1_+((0,1])$, 
respectively. These are homeomorphisms of $[0,1]$ that restrict to $C^1$ diffeomorphisms of $[0,1)$ and $(0,1]$, 
respectively, and to $C^{\ell-1}$ diffeomorphisms of $(0,1)$.

The vector fields $X_f$ and $Y_f$ do not necessarily coincide, and the Mather invariant captures this defect of coincidence. 
Given points $p,q$ in $(0,1)$, consider the maps $\psi_{X_f}^p : t\mapsto f_t(p)$ and $\psi_{Y_f}^q : t\mapsto f^t(q)$ from 
$\mathbb{R}$ to $(0,1)$. These are easily seen to be $C^{\ell}$ diffeomorphisms. Consider now the following change of coordinates:  
$$M^{p,q}_f := (\psi_{Y_f}^q)^{-1} \circ \psi_{X_f}^p : \mathbb{R} \to \mathbb{R}.$$ 
The fact that $f$ is the time-1 map of the flows of $X_f$ and $Y_f$ implies that $M^{p,q}_f$ 
commutes with the translation by 1. Hence, it induces a $C^{\ell}$ diffeomorphism of the circle $\R/\Z$, 
which we still denote by $M^{p,q}_f$. Note that if $p'=f_\alpha(p)$ (resp. $q'=f^\beta(q)$), then $\psi_{X_f}^{p'}=\psi_{X_f}^p\circ T_\alpha$  
(resp. $\psi_{Y_f}^{q'}=\psi_{Y_f}^{q}\circ T_\beta$).  

Thus, changing $p$ and $q$ amounts to pre and post composing $M^{p,q}_f$ with rotations. The class
$[ M_f^{p,q}] $ of $M^{p,q}_f$ 
modulo these $\mathrm{SO}(2,\mathbb{R})$-actions (on the left and right) is the Mather invariant of $f$, that we just 
denote by $M_f$. One says that the Mather invariant is trivial if $M_f$ coincides with the class of rotations. In view 
of the discussion above, this is equivalent to that $X_f$ and $Y_f$ coincide, which means that $f$ is 
{\em $C^1$-flowable}, that is, it arises as the time-1 map of the flow of a $C^1$ vector field on $[0,1]$. 

\medskip

More generally, one can interpret the Mather invariant of $f$ in terms of the $C^1$-centralizer of $f$; in particular, it provides information 
about the existence of a $k^{th}$ root of class $C^1$ of $f$ (that is, a $C^1$ diffeomorphism $g$ of $[0,1]$ such that $g^k=f$). Namely, 
if $\Gamma$ denotes the (closed) subgroup $\{t\in\R, f_t=f^t\}$ of $\R$ (which contains $\Z$), the centralizer $\mathrm{Z}_f^1$ of $f$ 
in $\Diff^1_+([0,1])$ is precisely $\{f_t,t\in\Gamma\}$. If, on the one hand, $\Gamma=\R$ (which means that $X_f$ and $Y_f$ coincide on $(0,1)$, 
or equivalently, that the Mather invariant of $f$ is trivial), then $f$ is $C^1$-flowable, and $\mathrm{Z}_f^1$ coincides with the corresponding flow. On the 
other hand, if $\Gamma\neq\R$ (which corresponds to the case of a nontrivial Mather invariant), then it is necessarily of the form $\frac1k\Z$ 
for some positive integer $k$, and $\mathrm{Z}_f^1$ is the infinite cyclic group generated by a $k^{th}$ root of $f$.  
In this case, it is not difficult to see that, given $p\in(0,1)$, the map $M_f^{p,p}$ commutes 
with the translation by $1/k$, and that $k$ is the largest integer such that $f$ has a $k^{th}$ root of class $C^1$ on $[0,1]$ (namely, $f^{1/k}=f_{1/k}$). 

\medskip

\begin{rem}
\label{r:mather-root}
The above implies that a $C^{\ell}$ diffeomorphism $f$ of $[0,1]$ without interior fixed point and with nontrivial Mather invariant cannot be conjugated 
to a nontrivial power $f^k$. Otherwise, it would have a $k^n$-th root of class $C^1$ for every positive integer $n$. But then, for any $p~\in~(0,1)$, the 
map $M_f^{p,p}$ would commute with the translation by $1/k^n$ for every $n$, hence it would be a translation, a contradiction. 
Another argument will be given in \S \ref{s:Vconj} using the relation between the Mather invariant 
(or $C^1$-flowability) and the asymptotic variation.
\end{rem}

\medskip

The Mather invariant was initially introduced to describe the $C^1$-conjugacy classes of $C^2$ diffeomorphisms 
of the interval without interior fixed points. Namely, two such diffeomorphisms are $C^1$-conjugated if and only if they have the 
same Mather invariant and their germs at the endpoints are $C^1$ conjugated. Though we will not use this precise result, we 
will use the ideas behind it, for instance in Proposition \ref{l:criterion} below.


\section{A proof of Theorem \ref{t:main}}
\label{s:distorted}

The proof of Theorem \ref{t:main} is based on the following proposition that gives 
a sufficient condition of $C^{\ell}$-conjugacy for a special type of diffeomorphisms.  
The proof, given in the next subsection, is an easy variation of ideas introduced in \S \ref{s:back}.  

\begin{prop}
\label{l:criterion}
Let $f$ and $g$ be $C^{\ell}$ diffeomorphisms of $[0,1]$ with the same finite set of fixed points, with $\ell \geq 1$. Suppose that their germs about 
each of these points are $C^{\ell}$-conjugated. Suppose also that $f$ embeds in the flow of a $C^{\ell}$ vector field, and the same holds for $g$. 
Then $f$ and $g$ are $C^{\ell}$-conjugated. 
\end{prop}

In view of this proposition, to prove Theorem \ref{t:main} it thus suffices to prove the density in $\Diff_P^1([0,1])$ of the set of diffeomorphisms $f$ 
such that $f$ and $g=f^k$ satisfy the hypotheses above for $\ell=\infty$ and every positive integer $k$.
The proof of this is divided in two steps:

\begin{enumerate}
\item The first step consists in approximating any $f \in \Diff_P^1([0,1])$ by $C^\infty$ diffeomorphisms with finitely 
many fixed points and with nice germs about these points, namely germs that are conjugated to all of their positive 
powers (these naturally arise inside the affine group -realized as a group of germs- and the ramified versions of it).

\item The second step consists in perturbing (in $C^1$ topology) every ``good'' 
element as above in order to make it ``$C^\infty$-flowable'' (meaning that 
it embeds in the flow of a $C^{\infty}$ vector field).
\end{enumerate}

Step 1 is divided in the three elementary lemmas below which are proved in \S \ref{s:elementary}, 
while Step~2 is the content of Proposition \ref{p:G4} and is proved in \S \ref{s:cancel}. The 
statements involve the following decreasing sequence of subsets of $G_0 := \Diff^1_P([0,1])$:

\begin{itemize}

\item $G_1$ is the subset of $\Diff^1_P([0,1])$ made of the diffeomorphisms that have only finitely many fixed points;

 \item $G_2$ is the subset of $G_1$ made of the diffeomorphisms that, about every $a\in \Fix(f)$, are conjugated 
 by the translation $T_a$ to $x\mapsto \frac{x}{1\pm x}$ or $x\mapsto\frac{x}{\sqrt{1\pm x^2}}$;

 \item $G_3$ is the subset of $C^\infty$ diffeomorphisms in $G_2$;

\item $G_4$ is the subset of diffeomorphisms of $G_3$ that are $C^\infty$-flowable.
 
\end{itemize}

Note that the germs involved in the definition of $G_2$ are conjugated to their positive powers by homotheties. Thus, for every $f\in G_4$ and every integer $k \geq 2$, 
the diffeomorphisms $f$ and $f^k$ satisfy the hypotheses of Proposition \ref{l:criterion}, and are thus $C^{\infty}$-conjugated. Therefore, Theorem \ref{t:main} is a direct 
consequence of the following facts:
 
\begin{lem}
\label{l:G1}
The set $G_1$ is dense in $G_0$ for the $C^1$ topology.
\end{lem}

\begin{lem}
\label{l:G2}
The set $G_2$ is dense in $G_1$ for the $C^1$ topology.
\end{lem}

\begin{lem}
\label{l:G3}
The set $G_3$ is dense in $G_2$ for the $C^1$ topology.
\end{lem}

\begin{prop}
\label{p:G4}
The set $G_4$ is dense in $G_3$ for the $C^1$ topology.
\end{prop}


\subsection{A proof of Proposition \ref{l:criterion}}
\label{s:cirterion}

This is an adaptation of Lemma 2.6 in \cite{EN21}. Let $f$ and $g$ be as in the statement of the proposition. In particular, there exist 
$C^{\ell}$ vector fields on $[0,1]$, denoted $X_f$ and $X_g$, whose time-$1$ maps are $f$ and $g$, respectively. According to \S \ref{s:back}, 
such vector fields are unique, and we denote by $(f^t)$ and $(g^t)$ their respective flows. Moreover, let us denote by $p_0=0 < p_1 < \ldots < p_n=1$ 
the common fixed points of $f$ and $g$. We will prove by induction on $k\in\{0,1,\ldots,n\}$ that $f$ and $g$ are $C^{\ell}$ conjugate on a neighborhood 
of $[p_0,p_k]$ in $[0,1]$ (by a map fixing $p_0,\dots,p_k$). 

The initial case $k=0$ follows directly from the hypothesis that the germs of $f$ and $g$ at $p_0$ are $C^{\ell}$-conjugate. 
Let us now assume that the result is true for some $k\in \{0,1,\ldots,n-1\}$. Namely, there exists a $C^{\ell}$ diffeomorphism $\phi_k$ between two 
neighborhoods of $[p_0,p_k]$ that conjugates $f$ to $g$ (and fixes $[p_0,p_k]$). Extend $\phi_k$ to a homeomorphism of $[p_0,p_{k+1}]$ 
in the unique possible way so that it still conjugates $f$ to $g$ on $[p_k,p_{k+1})$. Then the extension (still denoted $\phi_k$) is a $C^{\ell}$ 
diffeomorphism of $[p_0,p_{k+1})$. Now let $\hat\phi_k$ be a germ of $C^{\ell}$ diffeomorphism at $p_{k+1}$ conjugating $f$ to $g$, and 
similarly extend it to a $C^{\ell}$ diffeomorphism of $(p_k,p_{k+1}]$ still satisfying this property. By uniqueness of the generating vector 
fields on $[p_k,p_{k+1})$ (resp. $(p_k,p_{k+1}]$) we have $(\phi_k)_*X_f=X_g$ (resp. $(\hat\phi_k)_*X_f=X_g$) there, hence 
$(\hat\phi_k^{-1}\phi_k)_*X_f=X_f$ on $(p_k,p_{k+1})$, which means that the $C^{\ell}$ diffeomorphism $\hat\phi_k^{-1}\phi_k$ of $(p_k,p_{k+1})$ 
preserves the flow of $X_f$. This easily implies that $\hat\phi_k^{-1} \phi_k$ actually coincides with some time-$t_k$ map of this flow, so that 
$\phi_k=\hat\phi_k\circ f^{t_k}$ on $(p_k,p_{k+1})$. This allows extending $\phi_k$ to a $C^{\ell}$ diffeomorphism between neighborhoods 
of $[p_0,p_{k+1}]$ that conjugates $f$ to $g$, as required.


\subsection{Approximating by tamely behaved maps: proofs of Lemmas \ref{l:G1}, \ref{l:G2} and \ref{l:G3}}
\label{s:elementary}

\begin{proof}[Proof of Lemma \ref{l:G1}]
Let $f\in G_0$. Fix $\eps>0$, and consider the family $\mathcal{I}$ formed by the maximal intervals $I$ whose boundary points are fixed points 
of $f$ and where $\| (f-\id)|_I \|_1 \leq \eps/2$. On each interval $I$ of $\mathcal{I}$, replace $f$ by any diffeomorphism that is parabolic at the 
endpoints, has no fixed point in the interior and is $\eps/2$-close to the identity in $C^1$ norm. Let $\hat{f}$ be the new $C^1$ 
diffeomorphism thus obtained. By construction, $\Fix(\hat f)\subset \Fix(f)$. Moreover, $\hat{f}$ coincides with $f$ outside the 
intervals in $\mathcal{I}$, whereas on each $I \in \mathcal{I}$ we have
$$\| (\hat {f} - f)|_I\|_1 \leq \| (\hat{f} - \id)|_I \|_1 + \| ( \id - f )|_I  \|_1 \leq \frac{\varepsilon}{2} + \frac{\varepsilon}{2} = \varepsilon.$$
This shows that $\hat{f}$ is $\varepsilon$-close to $f$ in $C^1$ norm. 

Finally, we claim that $\Fix(\hat f)$ is finite. Indeed, assume by contradiction that this is not the case and let $p$ be an accumulation 
point of $\Fix(\hat f)$, say from the right. Then there exists an interval $[p,q]$ containing infinitely many fixed points of $\hat{f}$, hence 
of $f$. By slightly approaching $q$ to $p$ if necessary, we may assume that $\| (f-\id)|_{[p,q]} \|_1 \leq \varepsilon / 2$. However, this 
is in contradiction with the definition of $\hat{f}$.
\end{proof}

\begin{proof}[Proof of Lemma \ref{l:G2}] We first fix some notation. 
We call $q_1$ and $q_2$ the local diffeomorphisms at the origin defined by
$$q_1 (x) := \frac{x}{1-x}  \qquad \mbox{ and } \qquad q_2 (x) := \frac{x}{\sqrt{1-x^2}}.$$ 
Note that 
$$q_1 ^{-1}(x) = \frac{x}{1+x}  \qquad \mbox{ and } \qquad q_2^{-1} (x) = \frac{x}{\sqrt{1+x^2}}.$$ 
Also, given $\varepsilon > 0$, let $\rho = \rho_{\varepsilon}: [-\varepsilon,\varepsilon] \to [0,1]$  be a  $C^{\infty}$ 
function such that $\rho|_{[-\varepsilon/4,\varepsilon/4]} \equiv 1$, $\rho|_{[-\varepsilon,-3\varepsilon/4]} \equiv 0$, 
$\rho|_{[3\varepsilon/4,\varepsilon]} \equiv 0$ and $\| D\rho \|_{\infty} \leq \frac{3}{\varepsilon}$. 

\begin{center}
	\includegraphics[scale=0.31]{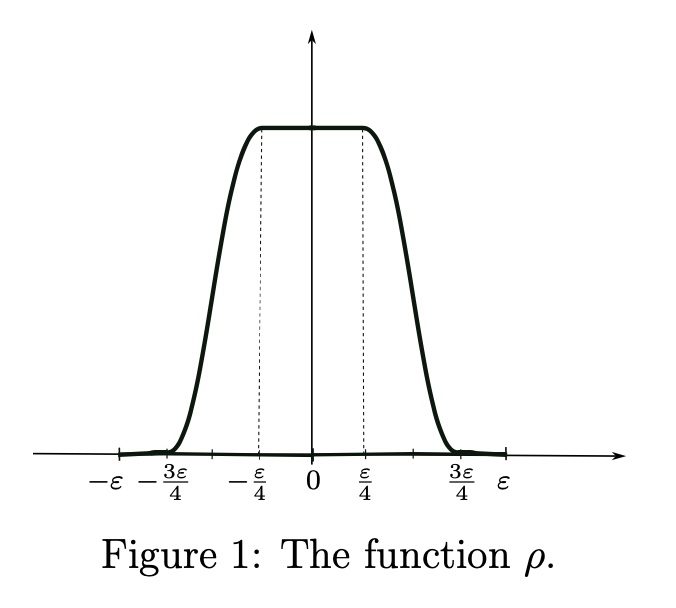}\\
\end{center}

Let now $f \in G_1$, and let $\varepsilon>0$ be small enough so that the $\varepsilon$-neighborhoods of the fixed points 
of $f$ are disjoint. Given an interior fixed point of $f$, say that $a$ is of type 1 (resp. type~2) if the graph of $f$ 
does not cross (resp. crosses) the diagonal at $a$. Let $q_{1,a}$ and $q_{2,a}$ be the conjugates of $q_1$ and 
$q_2$ by the translation $T_a$. Let us define $\tilde{f}$ on $[a-\varepsilon,a+\varepsilon]$ as follows:
 
\vspace{0.2cm}

\noindent - If $a$ is of type 1 then, on $[a-\varepsilon,a+\varepsilon]$, let $\hat{f}$ be equal to $\rho \, q_{1,a} + (1-\rho) f$ 
if the graph of $f$ is above the diagonal about $a$, and let $\hat{f} := \rho \, q_{1,a}^{-1} + (1-\rho) f$ otherwise.
 
 \vspace{0.2cm}
 
\noindent - If $a$ is of type 2 then, on $[a-\varepsilon,a+\varepsilon]$, let $\hat{f}$ be equal to $\rho \, q_{2,a} + (1-\rho) f$ if the graph 
of $f$ is below the diagonal on the left of $a$ (and above on the right), and let $\hat{f} := \rho \, q_{2,a}^{-1} + (1-\rho) f$ otherwise.

\vspace{0.2cm}

Perform this modification about each interior fixed point of $f$, and proceed analogously with each endpoint of $[0,1]$ on the 
corresponding one-sided neighborhood (on each of these, one can use either $q_1$ or $q_2$ if the graph is above the diagonal, and 
the inverses otherwise). Using that $0 \leq \rho \leq 1$, one can then easily check that this procedure does not introduce new fixed points. 
Therefore, letting $\hat{f}$ be equal to $f$ outside the $\varepsilon$-neighborhoods of the fixed points of $f$, we get a $C^1$ diffeomorphism 
with the desired prescribed germs at the fixed points. In other words, $\hat{f}$ belongs to $G_2$. Note that $\hat{f}$ depends on $\varepsilon$.  

Finally, we claim that if $\varepsilon$ goes to zero, then the corresponding elements $\hat{f}$ constructed above converge to $f$ 
in $C^1$ topology. This is to be checked on the $\varepsilon$-neighborhoods of the fixed points $a$ of $f$ 
(since $\hat{f}$ coincides with $f$ outside). Letting $q$ be equal to 
$q_{1,a}^{\pm 1}$ or $q_{2,a}^{\pm 1}$ according to the case, we compute:
$$D \hat{f} - Df = \rho (Dq - Df) + (q-f) D \rho.$$
Since both $D(q-f)$ is of the form $o(1)$ about $a$ and $\rho \leq 1$, the first term above $\rho \, (Dq - Df)$ goes to zero as 
$\varepsilon \to 0$.  Also, $q-f$ is of the form $o(x-a)$ about $a$, and therefore 
$$| (q-f) D \rho | \leq \frac{3}{\varepsilon} |q-f| = o(1),$$
which shows the desired convergence to $0$.
\end{proof}

\begin{proof}[Proof of Lemma \ref{l:G3}] We will show that every element $f \in G_2$ with no fixed point in $(0,1)$ can be 
approximated by elements in $G_3$ with no fixed points in $(0,1)$. To deal with a general $f \in G_2$, it suffices to apply the 
same method on the closure of each of the (finitely many) connected components of $[0,1] \setminus \mathrm{Fix} (f)$. 

Let hence $f \in G_2$ have no fixed point in $(0,1)$. Fix $\varepsilon > 0$ small enough so that $f$ coincides with one of $q_1^{\pm 1}$ 
or $q_2^{\pm 1}$ on $[0,2\varepsilon]$ and with a conjugate of one of them on $[1-2\varepsilon,1]$. 
The function $Df$ restricted to $[\varepsilon,1-\varepsilon]$ 
can be approximated by a sequence of $C^{\infty}$ functions $\varphi_n$ (e.g. use Berstein's polynomials). Let $\hat{f}_n$ be the map 
obtained by integrating $\varphi_n$ on $ [ \varepsilon , 1 - \varepsilon ] $ with the same value as $f$ at $ \varepsilon$. More precisely,  
$$\hat{f}_n (x) := f( \varepsilon) + \int_{ \varepsilon}^x \varphi_n (x) dx, \qquad x \in [\varepsilon,1-\varepsilon].$$
Note that $\hat{f}_n$ converges to $f$ in the $C^1$ sense on $[\varepsilon,1-\varepsilon]$, meaning that both 
$$\sup_{x \in [\varepsilon,1-\varepsilon]} |\hat{f}_n (x) - f (x)| 
\qquad \mbox{ and } \qquad 
\sup_{x \in [\varepsilon,1-\varepsilon]} |D \hat{f}_n (x) - D f(x)|$$
converge to zero as $n$ goes to infinity.
 
Finally, let $f_n$ be defined as the interpolation along $\rho$ between $f$ and $\hat{f}_n$ on the intervals 
$[\varepsilon,2\varepsilon]$ and $[1-2\varepsilon,1-\varepsilon]$ (cf. Figure \ref{f:interp}). More precisely, 
$$f_n(x)
=
\begin{cases} 
      f(x) & \mbox{if } x \in [0,\varepsilon] , \\
            \rho(x-\varepsilon) f(x) + (1-\rho(x-\varepsilon)) \hat{f}_n(x) & \mbox{if } x \in [\varepsilon,2\varepsilon], \\
      \hat{f}_n(x) & \mbox{if } x \in [2 \varepsilon,1-2\varepsilon], \\
                (1 -  \rho(x-1+2\varepsilon)) f(x) + \rho(x-1+2\varepsilon) \hat{f}_n(x) & \mbox{if } x \in [1-2\varepsilon,1-\varepsilon], \\
            f(x) & \mbox{if } x \in  [1-\varepsilon,1]. \\
\end{cases}
$$

\begin{center}
	\includegraphics[scale=0.34]{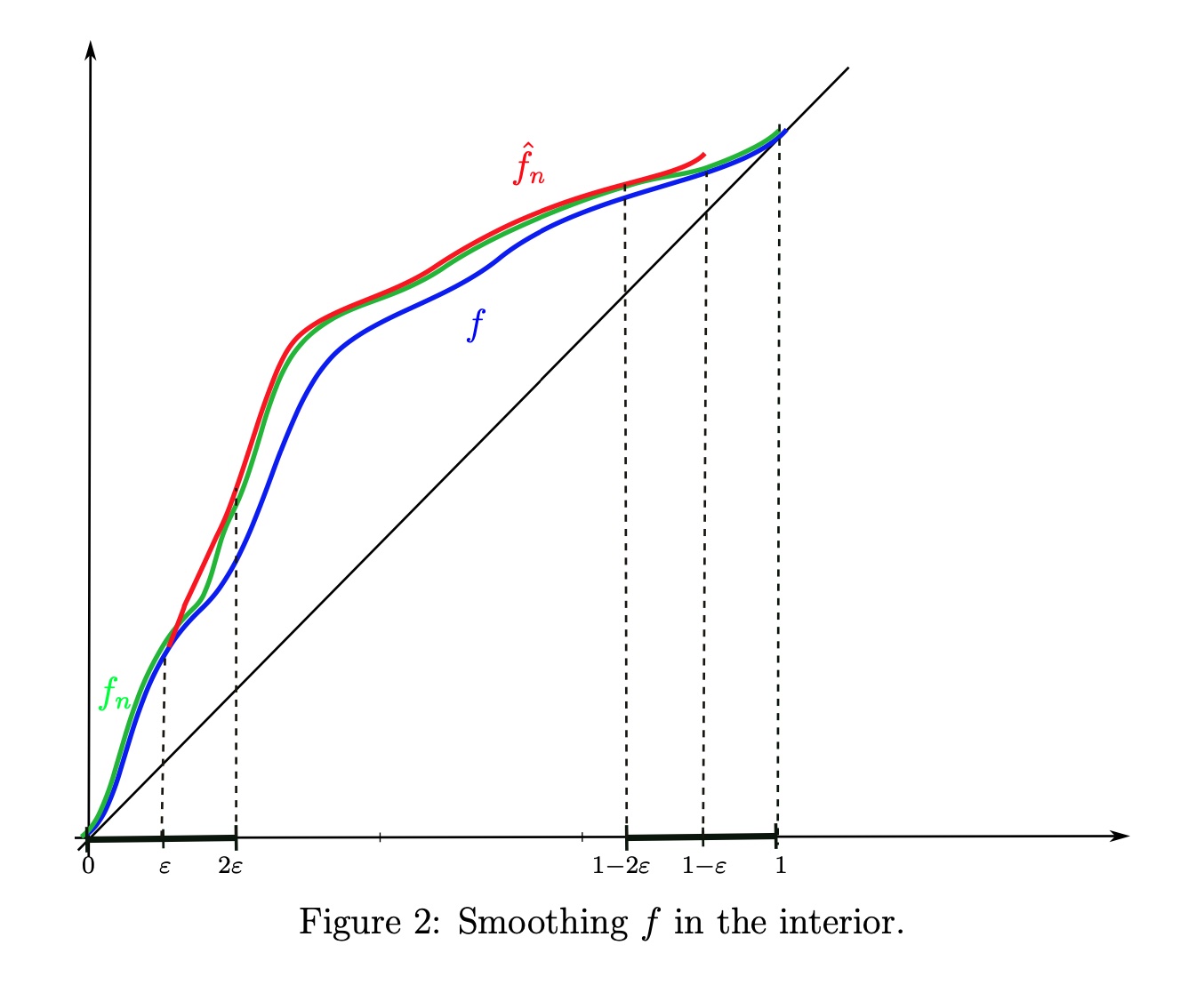}\\
\end{center}

This is a $C^{\infty}$ map which, actually, belongs to $G_3$ for large-enough $n$. We claim that $f_n$ converges to $f$ in $C^1$ topology. 
This is obvious outside $[\varepsilon,2\varepsilon] \cup [1-2\varepsilon,1-\varepsilon]$, since on $[0,\varepsilon] \cup [1-\varepsilon,1]$  
one has $f_n = f$ for all $n$, whereas on $[2\varepsilon,1-2\varepsilon]$ we have $f_n = \hat{f}_n$, and the latter converges to $f$ 
in the $C^1$ sense. 

Let us now deal with $[\varepsilon,2\varepsilon]$, 
the case of $ [1-2\varepsilon,1-\varepsilon]$ being analogous. For $x \in [\varepsilon,2\varepsilon]$, we compute: 
$$Df_n(x) - D \hat{f}_n (x) = (f - \hat{f}_n)(x) \, D\rho (x-\varepsilon) + (Df - D\hat{f}_n)(x) \, \rho (x-\varepsilon).$$
Since $\varepsilon$ (hence $\rho$) is fixed and $\hat{f}_n$ converges to $f$ in the $C^1$ sense on $[\varepsilon,1-\varepsilon]$, the two 
terms in the sum above converge to zero, as desired.
\end{proof}

\subsection{Cancelling the Mather invariant: a proof of Proposition \ref{p:G4}}
\label{s:cancel}

This is by far the hardest part of the proof of Theorem \ref{t:main}, and crucially involves the notions introduced 
in \S \ref{s:back} such as generating vector fields and the Mather invariant. Each $f \in G_3$ is 
a $C^\infty$ diffeomorphism with finitely many fixed points near which it is the time-1 map of a $C^\infty$ vector field. 
We will perform a $C^1$-small (yet $C^\infty$) perturbation of $f$ supported on compact subintervals of $(0,1)$ 
arbitrarily close to the fixed points but not containing them in such a way that the resulting diffeomorphism is 
$C^1$-flowable, or, in other words, has a trivial Mather invariant. This perturbed diffeomorphism 
will then be automatically $C^\infty$-flowable. Indeed, this is because of the uniqueness of the $C^1$ generating vector field 
(whenever it exists) and the fact that in our context, it is $C^\infty$ near the fixed points (since we have not modified 
the starting diffeomorphism therein), and $C^\infty$ away from them since it is invariant under its time-1 map which 
is assumed to be $C^\infty$. \medskip

The claim that this perturbation is possible is similar to \cite[Proposition 1.4]{BCVW08} (applied to each component of~$f$). 
The only difference is that our $f$ is parabolic at the endpoints, whereas \cite{BCVW08} deals with hyperbolic (and even affine) 
maps at the endpoints. However, we will see that this does not play such an important role (concretely, the difference appears in 
the proof of Lemma \ref{l:tricky} below). 

\medskip

Let us now enter into the details. Fix $f\in G_3$. For simplicity, we will assume that $f$ has no interior fixed point (since we will perform perturbations 
away from the fixed points, which are in finite number, so everything matches up smoothly). To fix ideas, assume that $f(x)-x>0$ for all $x$ in $(0,1)$, 
and  denote by $X_f$ and $Y_f$ its left and right generating vector fields, respectively.

The argument works as follows. In Lemma \ref{l:modify} and its Corollary \ref{c:modify} below, we will show that modifying $f$ on finitely many disjoint 
fundamental domains closer and closer to $0$ we can change
(a prescribed lift of) its Mather invariant, the resulting one being a composition of the initial one 
with (lifts of) prescribed circle diffeomorphisms of a special type (namely, equal to the identity on some interval of the circle). Using a fragmentation lemma 
for circle diffeomorphisms (cf. Lemma~\ref{l:fragmentation} below), this will actually allow cancelling the Mather invariant of $f$. The tricky part will consist in 
showing that if the fragments are $C^1$-small and the perturbation occurs close enough to~$0$, then the modified 
diffeomorphism is $C^1$-close to $f$. This is more or less the content of Lemmas \ref{l:key} and \ref{l:tricky} below.

\medskip

From now on, we fix some $p\in(0,1)$ and we define the diffeomorphisms $\psi_{X_f}:=\psi_{X_f}^p$ and $\psi_{Y_f}:=\psi_{Y_f}^p$ from $\R$ to $(0,1)$ 
as in \S \ref{s:back} (so that the images of $0$ under these two maps are~$p$). Given a diffeomorphism $\phi$ of $\R$ supported in an interval 
$(\alpha-1,\alpha)$, we denote by $\tilde \phi$ the unique diffeomorphism of $\R$ coinciding with $\phi$ on $[\alpha-1,\alpha]$ and commuting with the unit translation.

\begin{lem}
\label{l:modify} Let $\phi$ be a diffeomorphism of $\R$ supported in some open interval of the form $(\alpha-1,\alpha)$, with $\alpha\le 0$, 
and  let $h := \psi_{Y_f} \phi \, \psi_{Y_f}^{-1}$. (This is a diffeomorphism of $[0,1]$ supported in the fundamental domain $[f^{-1}(q),q]$ of $f$, where 
$q:=\psi_{Y_f}(\alpha)=f^\alpha(p)$.) If $g=f\circ h$, then $\psi_{Y_g}^p = \psi_{Y_f}^p$ on $[\alpha,+\infty)$ and $M^{p,p}_g\circ T_\tau= \tilde\phi\circ M_f^{p,p}$ for some $\tau\in\R$.
\end{lem}

\begin{proof} This claim is similar to \cite[Lemma 1.5]{BCVW08} (but with different notations) and to \cite[Remark 2.3]{EN21}. 
We adapt the argument of the latter reference for self-containedness. 

\begin{center}
	\includegraphics[scale=0.27]{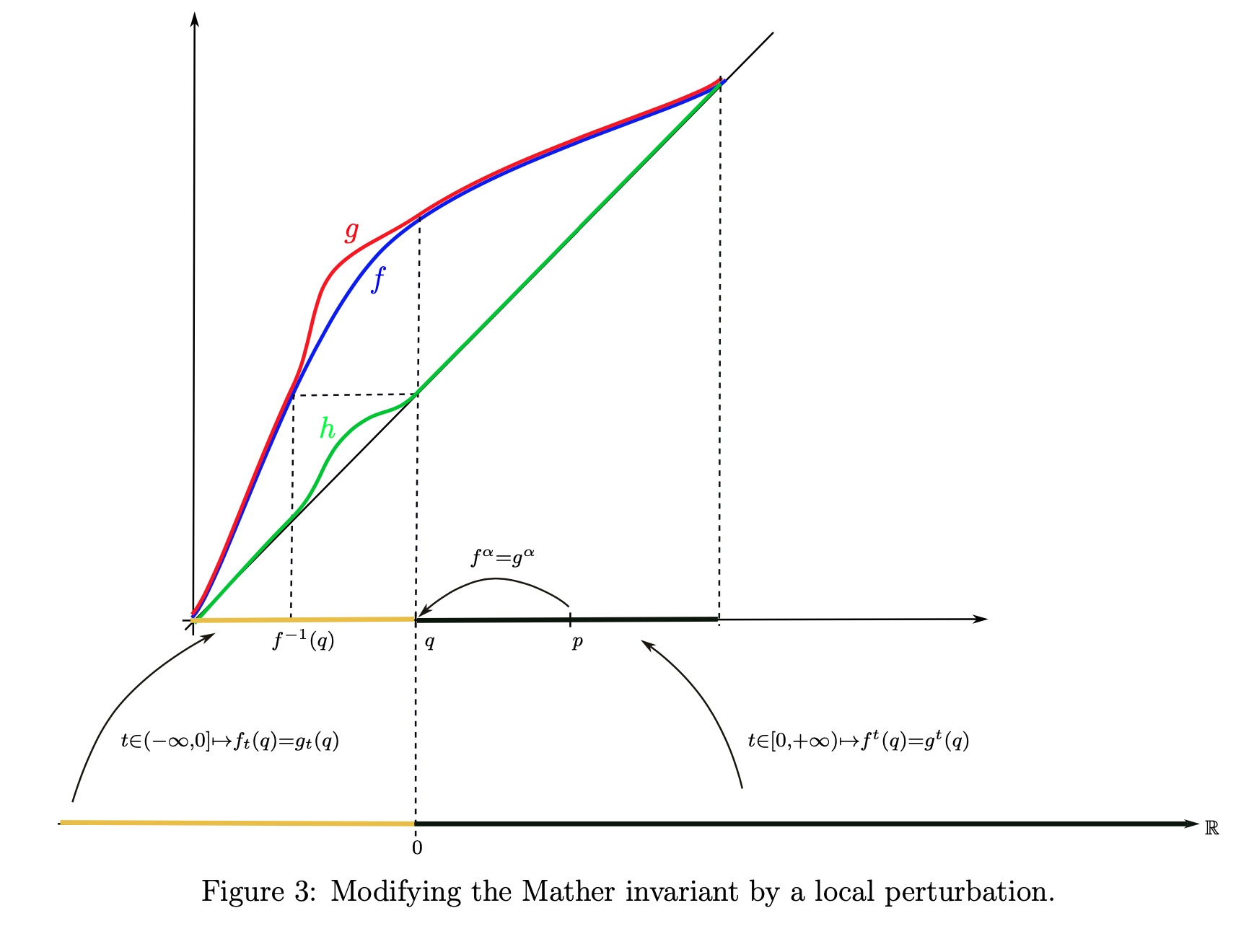}\\
\end{center}

Since $g$ (resp. $g^{-1}$) coincides with $f$ on $[0,f^{-1}(q)]$ (resp. with $f^{-1}$ on $[f(q),1]$), 
necessarily $X_g$ (resp. $Y_g$) coincides with $X_f$ (resp.  $Y_f$) on $[0,q]$ (resp. $[q,1]$). This easily implies that: 
\begin{equation}
\label{e:fg1}
(\psi_{X_f}^q)^{-1} \text{ and } (\psi_{X_g}^q)^{-1} \text{  define the same diffeomorphism from } (0,q] \text{ to } (-\infty,0], 
\end{equation}
\begin{equation}
\label{e:fg2}
(\psi^{q}_{Y_f})^{-1} \text{ and } (\psi^{q}_{Y_g})^{-1} \text{ define the same diffeomorphism from } [q,1) \text{ to } [0,+\infty).
\end{equation}

In particular, since $-\alpha$ is nonnegative, 
$$p=f^{-\alpha}(q)=\psi_{Y_f}^q(-\alpha) = \psi_{Y_g}^q(-\alpha)=g^{-\alpha}(q),$$
hence
$g^\alpha(p)=f^\alpha(p)=q$. 
Therefore, for all $t \in \R$, one has 
\begin{equation}
\label{e:Ta}
\psi_{Y_f}^q(t)=f^t(q)=f^t(f^\alpha(p))=\psi_{Y_f}^p(t+\alpha)=\psi^p_{Y_f}\circ T_\alpha (t), \end{equation}
and the same holds replacing $f$ by $g$:
\begin{equation}
\label{e:Tabis}
\psi_{Y_g}^q(t)=g^t(q)=g^t(g^\alpha(p))=\psi_{Y_g}^p(t+\alpha)=\psi^p_{Y_g}\circ T_\alpha (t). \end{equation}

Similarly, 
\begin{equation}
\label{e:Tb}
\psi_{X_f}^{q}=\psi_{X_f} \circ T_\beta \quad\text{ and } \quad\psi_{X_g}^{q}=\psi_{X_g} \circ T_\gamma 
\end{equation} 
for some nonpositive real numbers $\beta,\gamma$ 
such that $f_\beta(p)=q$ and $g_\gamma(p)=q$ (note that $\beta$ and $\gamma$ do not necessarily coincide, 
because $X_f$ and $X_g$ differ on $[q,p]$).

\medskip

Recall finally that: 

\medskip

\noindent $(*)$ For every $r$, the maps $\psi^{r}_{X_f}$ and $\psi^{r}_{Y_f}$ (resp. $\psi^{r}_{X_g}$ and $\psi^{r}_{Y_g}$) 
conjugate $f$ (resp. $g$) to the unit translation $T_1$.  
\smallskip

Putting all of this together, on the interval
$$(M^{p,p}_f)^{-1}([\alpha-1,\alpha])
=\psi_{X_f}^{-1}\psi_{Y_f}([\alpha-1,\alpha])
= \psi_{X_f}^{-1} ([f^{\alpha-1}(p),f^\alpha(p)])
=\psi_{X_f}^{-1}([f^{-1}(q),q]) 
= [\beta-1,\beta],$$
we obtain 
\begin{align*}\label{e:comput}
\phi\circ M^{p,p}_f&= \psi_{Y_f}^{-1}\circ h \circ \psi_{Y_f}\circ M^{p,p}_f\\
&=\psi_{Y_f}^{-1}\circ h \circ \psi_{X_f} \hspace{0.15cm}\quad\qquad \qquad \qquad \qquad \qquad \qquad\text{by definition of $M_f^{p,p}$,} \\
&= \psi_{Y_f}^{-1}\circ (f^{-1}\circ g) \circ \psi_{X_f} \hspace{0.04cm} \quad\qquad \quad \qquad \qquad \qquad\text{by definition of $g$,}\\
&= T_{-1} \circ \psi_{Y_f}^{-1}\circ g\circ \psi_{X_f} \hspace{0.4cm}\qquad \qquad \qquad \qquad \qquad\text{according to $(*)$,}\\
&= T_{-1} \circ T_{-\alpha}\circ (\psi_{Y_f}^{q})^{-1}\circ g\circ \psi_{X_f}^q\circ T_{-\beta}\hspace{0.06cm}\qquad\qquad\text{according to \eqref{e:Ta} and \eqref{e:Tb}}\\
&=T_{-1} \circ T_{-\alpha}\circ (\psi_{Y_g}^{q})^{-1}\circ g\circ \psi_{X_g}^q\circ T_{-\beta}\hspace{0.1cm}\qquad\qquad\text{according to \eqref{e:fg1} and \eqref{e:fg2}}\\
&= T_{-1} \circ T_{-\alpha}\circ T_1\circ (\psi_{Y_g}^{q})^{-1} \psi_{X_g}^q\circ T_{-\beta}\hspace{0.32cm}\qquad\qquad\text{according to $(*)$,}\\
&=T_{-\alpha}\circ( \psi_{Y_g}^{q})^{-1}\circ \psi_{X_g}^q \circ T_{-\beta}\\
&=T_{-\alpha}\circ T_\alpha\circ (\psi_{Y_g}^p)^{-1} \circ \psi_{X_g}^p \circ T_\gamma\circ T_{-\beta}\hspace{0.11cm}\qquad\qquad\text{according to \eqref{e:Tabis} and \eqref{e:Tb}}\\
&= M^{p,p}_g\circ T_{\gamma-\beta},
\end{align*}
which gives the announced result. 
\end{proof}

By induction, one obtains the following analog of \cite[Corollary 1.8]{BCVW08}:

\begin{cor}
\label{c:modify}
Consider a finite sequence $(\alpha_i)$ of real numbers, where $i\in\{1,\dots,l\}$, such that $\alpha_{i+1}<\alpha_i-1$ and $\alpha_1\le0$. 
For every $i$, let $\phi_i$ be a diffeomorphism of $\R$ supported in $(\alpha_i-1,\alpha_i)$, and let $h_i=\psi_{Y_f} \phi_i\psi_{Y_f}^{-1}$. 
(These are diffeomorphisms of $[0,1]$ supported in disjoint fundamental domains of $f$ closer and closer to $0$.)  If 
$g=f\circ h_1\circ h_2\circ...\circ h_l = f\circ h_l\circ\dots\circ h_1$, then 
$M_g^{p,p}\circ T_\tau = \tilde\phi_1\circ \dots\circ \tilde\phi_l \circ M_f^{p,p}$ for some $\tau\in\R$.
\end{cor}

\begin{proof} We write the proof only to emphasize the importance of the fact that the maps $h_1,\dots,h_l$ commute 
(which is ensured by the disjointness of their support, 
which will be also used in the proof of Lemma \ref{l:key} below). 
For every $i\in\{0,\dots,l-1\}$, define $g_i := f\circ h_l\circ\dots\circ h_{l-i}$. 
Let us prove by induction that, for every $i\in\{0,\dots,l-1\}$, it holds 
$M_g^{p,p}\circ T_{\tau_i} = \tilde\phi_{l-i}\circ \dots\circ \tilde\phi_l \circ M_f^{p,p}$ 
for some $\tau_i\in\R$ and $\psi_{Y_{g_i}}=\psi_{Y_f}$ on $[\alpha_i,+\infty)$.
The first step ($i=0$) is precisely the content of Lemma~\ref{l:modify}. Now, 
if we assume that the statement is true for some $i\in\{0,\dots,l-2\}$, 
then we apply Lemma~\ref{l:modify} replacing $f$ by $g_i$ and $\phi$ by $\phi_{l-i-1}$. 
It is crucial here to observe that 
$\psi_{Y_{g_i}}\phi_{l-i-1}\psi_{Y_{g_i}}^{-1}$ is indeed the diffeomorphism 
$\psi_{Y_f} \phi_{l-i-1} \psi_{Y_f}^{-1} = h_i$ 
due to the second part of the induction hypothesis.
\end{proof} 

Now we need to control the $C^1$ size of the perturbation obtained by the process above.

\begin{lem}
\label{l:key}
For every neighborhood $U$ of $f$ in $\Diff^1_+ ([0,1])$
there exists $\eps>0$ and $\beta \leq 0$ with the following property: Given real numbers $\alpha_i$ and 
diffeomorphisms $\phi_i$ as in the previous corollary 
such that $\phi_i$ is $\eps$-close to the identity in $C^1$ norm and $\alpha_1\le \beta$,
the diffeomorphism $f\circ h_1\circ\dots \circ h_l$ 
belongs to $U$ (where $h_i$ is defined as above).
\end{lem}

This is analogous to \cite[Lemma 1.8]{BCVW08}. We fix a neighborhood $U_0$ of the identity map of $[0,1]$ such that, 
if $h_1,\dots,h_l$ are disjointly supported and belong to $U_0$, then $f\circ h_1\circ \dots \circ h_l$ belongs to $U$. 
The claim is then a direct consequence of the following key lemma.

\begin{lem}
\label{l:tricky}
For every neighborhood $U_0$ of the identity in $\Diff^1_+ ([0,1])$, 
there exist $\eps>0$ and  $\beta \leq 0$ such that, if $\alpha\le \beta$, 
then for every diffeomorphism $\phi$ of $\R$ supported in $(\alpha-1,\alpha)$ 
that is $\eps$-close to the identity in $C^1$ norm,
the diffeomorphism $h := \psi_{Y_f} \phi\psi_{Y_f}^{-1}$ of $[0,1]$ belongs to $U_0$.
\end{lem}

\begin{proof}
As we already mentioned, although this is analogous to \cite[Lemma 1.9]{BCVW08}, our proof is different since, in our setting, 
the map $f$ we are dealing with is not affine at the endpoints.

Let $\delta>0$ be such that $U_0$ contains all diffeomorphisms  
$g$ such that $\|\log Dg\|_\infty<\delta$. Fix $\beta \leq 0$ such that, if $q=\psi_{Y_f}(\beta)$, then 
$V (f;[0,q])  < \delta / 2$. 
Let $\alpha\le \beta$. Consider $n\in\N$ such that $\alpha+n\in[\beta-1,\beta]$. Given a diffeomorphism $\phi$ of $\R$ supported in $(\alpha-1,\alpha)$, one has 
$$h=\psi_{Y_f} \phi\psi_{Y_f}^{-1}=f^{-n}f^n\psi_{Y_f} \phi\psi_{Y_f}^{-1}f^{-n}f^n=f^{-n}\psi_{Y_f}T_n\phi T_{-n}\psi_{Y_f}^{-1}f^n = f^{-n} h_n f^n$$
 where $h_n$ is defined as $ \psi_{Y_f}T_n\phi T_{-n}\psi_{Y_f}^{-1}$. Note that the latter  map is supported in 
$$\psi_{Y_f}((\alpha+n-1,\alpha+n))\subset\psi_{Y_f}([\beta-2,\beta])= [f^{-2}(q),q].$$ 

To estimate $\log(Dh)$, we proceed as follows. By the chain rule,
\begin{align*}
\log Dh = \log D h_n \circ f^n + \log\frac{Df^n}{Df^n \circ h} 
\end{align*}
hence
\begin{equation}\label{eq:estimate}
\|\log Dh\|_\infty\le \|\log D h_n \|_\infty+ \left\|\log\frac{Df^n}{Df^n \circ h} \right\|_\infty.
\end{equation}

The second term above is treated via a classical control of distortion argument. Namely, for every $x$ in the support of $h$, the points 
$x$ and $h(x)$ are contained in the same fundamental domain of~$f$. As a consequence, when $i$ varies from $0$ to $n-1$, the intervals 
with endpoints $f^i(x), f^i(h(x))$ are disjoint and contained in $[0,f^{n+\alpha}(p)]\subset [0,q]$, and therefore
\begin{align}
\left|\log\frac{Df^n}{Df^n \circ h}(x)\right| &=\left| \log Df^n(x) - \log Df^n (h(x))\right|\notag\\
&=\left|\sum_{i=0}^{n-1} \log Df (f^i(x))-\log Df(f^i \circ h(x))\right|
\le V (f; [0,q]) < \frac{\delta}{2}.
\label{eq:estimate2}
\end{align}

To deal with the first term in (\ref{eq:estimate}), first note that $h_n = \psi_{Y_f}  \phi_n \psi_{Y_f}^{-1}$, where 
$\phi_n := T_n\phi T_{-n}$. Denote $\omega$ the modulus of continuity of $\log (D \psi_{Y_f})$ on $[\beta-2,\beta]$, that is,   
$$\omega (\kappa) := \sup \{ |\log D\psi_{Y_f} (y) - \log D \psi_{Y_f} (z) | : \hspace{0.1cm} y,z \in [\beta-2,\beta],  \hspace{0.1cm} |y - z | \leq \kappa \}.$$
Note that $\omega (\kappa)$ goes to 0 as $\kappa$ goes to 0. By the chain rule, for $x \in [f^{-2}(q),q]$, 
\begin{eqnarray*}
Dh_n (x) 
&=& 
\frac{D \psi_{Y_f} (\phi_n \psi_{Y_f}^{-1} (x))}{D \psi_{Y_f} ( \psi_{Y_f}^{-1} (x))} \cdot D \phi_n (\psi_{Y_f}^{-1} (x)).
\end{eqnarray*}
Together with the definition of $\phi_n$, this immediately gives 
$$\| \log Dh_n  \|_{\infty} 
\leq 
\omega (\| \phi - id\|_{\infty}) + \| \log D \phi \|_{\infty}.
$$
Now choose $\varepsilon$ small enough so that, if $\phi$ is $\eps$-close to the identity in $C^1$ norm, then both terms above are smaller 
than $\delta / 4$. Then $\| \log Dh_n  \|_{\infty} < \delta /2$, which together with (\ref{eq:estimate}) and (\ref{eq:estimate2}) gives 
$\| \log Dh_n  (x)\|_{\infty} < \delta$. By the choice of $\delta$, this implies that $h$ belongs to $U_0$, as claimed.
\end{proof}

\medskip

To conclude the proof of the density of $G_4$ in $G_3$ (\emph{i.e.} of Proposition \ref{p:G4}), we need the 
following corollary of a classical (and much more general) Fragmentation Theorem (see for instance \cite{Ba97}):

\medskip

\begin{lem}
\label{l:fragmentation}
For each neighborhood $W$ of the identity in $\Diff_+^1(\mathbb{S}^1)$ and every $C^\infty$ diffeomorphism $\varphi$ of the circle, 
there exists a positive integer~$l$ and diffeomorphisms $\varphi_1,\dots,\varphi_l$ in $W$ that are equal to the identity on 
some interval and satisfy $\varphi = \varphi_1\circ\dots\circ\varphi_l$. 
\end{lem}

\bigskip

Consider now our initial $f \in G_3$, a $C^1$-neighborhood $U$ of $f$ and the corresponding $\eps>0$ and $\beta\in\R$ given by Lemma \ref{l:key}. 
Let $W$ be the $C^1$-neighborhood of the identity of the circle formed by the diffeomorphisms that are $\eps$-close to the identity in $C^1$ norm. 
Use Lemma \ref{l:fragmentation} to write the circle diffeomorphism (induced by) $(M_f^{p,p})^{-1}$ as a composition $\varphi_1\circ \dots\circ \varphi_l$ 
with each $\varphi_i\in W$ and equal to the identity about some $\bar \alpha_i \in \mathbb{S}^1$. 
Choose lifts $\alpha_i$ of the points $\bar \alpha_i$ such that $\alpha_{i+1}<\alpha_i-1$ for every~$i$ and $\alpha_1\le \beta$. 
For each $i$ define the unique diffeomorphism $\phi_i$ of $\R$ supported on $[\alpha_i-1,\alpha_i]$ and induced by (a lift of)
$\varphi_i$ therein, and extend it to a diffeomorphism $\tilde{\phi}_i$ of $\R$ that commutes with the unit translation. 
Finally, let $h_i := \psi_{Y_f}\circ \phi_i\circ\psi_{Y_f}^{-1}$. According to Lemma~\ref{l:key}, the diffeomorphism $g:=f\circ h_1\circ \dots\circ h_l$ belongs to $U$, 
and according to Lemma~\ref{c:modify}, one has the equality $M_g^{p,p} \circ T_\tau = \tilde\phi_1\circ \dots\circ \tilde\phi_l\circ M_f^{p,p}$ for some $\tau\in\R$. 
This implies that $M_g = [ \varphi_1\circ\dots\circ\varphi_l \circ M_f^{p,p}] = [ (M_f^{p,p})^{-1} \circ M_f^{p,p}]$ is trivial, so $f$ is $C^1$-flowable. Since $g$ 
coincides with $f$ near the endpoints, this implies that $g$ is actually $C^\infty$-flowable, which concludes the proof of Proposition \ref{p:G4}.
 

\section{A first proof of Theorem \ref{t:generic}}
\label{s:generic-short}

The first proof we give of Theorem \ref{t:generic} is somewhat deceitful: it is based on a black box coming from \cite{CC02} and \cite{BMNR17} 
which is summarized in the next two results.

\begin{thm} {\em Let $f,h$ be $C^1$ diffeomorphisms of $[0,1]$ such that $hfh^{-1} = f^k$ for some $k \geq 2$. If $I$ is a component 
of $f$, then $h$ fixes $I$.}
\label{t:CC}
\end{thm}

\begin{thm} {\em Let $f,h$ be $C^1$ diffeomorphisms of $[0,1]$ such that $hfh^{-1} = f^k$ for some $k \geq 2$. If $f$ 
has no fixed point on $(0,1)$, then the action of the group $\langle f,h \rangle$ is topologically conjugated to that of a subgroup 
of the affine group. In particular, $h$ has a unique fixed point $p$ at the interior. Moreover, at this point, one has $Dh (p) = k$.}
\label{t:BMNR}
\end{thm}

 Although stated in a different terminology, the first theorem above is essentially proved in \cite{CC02} for $k=2$, yet the proof works verbatim for any 
 $k \geq 2$; see also \cite{Na10,GL11} as well as \cite[Proposition 1.5]{BMNR17}. The second theorem, which is strongly inspired by \cite{CC02}, is 
 a direct consequence of some of the the main results from \cite{BMNR17} (namely Proposition 1.6 and Theorem 1.7 therein).  We will only use 
 the one concerning the rigidity of the multiplier of $h$ at the fixed point, which is highly nontrivial since, in general, the topological 
 conjugacy to the affine action in the statement fails to be differentiable (see \cite[Proposition 4.15]{BMNR17} in this regard).
 
In order to apply these results, we will need the following elementary lemma.

\begin{lem} 
The subset $\mathcal{S}$ of $\Diff^1_P ([0,1])$ made of the diffeomorphisms $f$ for which 
the function $x \mapsto f(x)-x$ changes its sign in every neighborhood 
of the origin is residual in $\Diff^1_P ([0,1])$.
\end{lem}

\begin{proof} 

For each $n \in \N$, consider the subset $\mathcal{S}_n$ made of the diffeomorphisms 
$f \in \Diff^1_P ([0,1])$ for which $x \mapsto f(x)-x$ has a change of sign in $(0,1/n)$. Obviously, 
$\mathcal{S} = \bigcap_{n \in \N} \mathcal{S}_n$, so it suffices to show that each $\mathcal{S}_n$ is 
open and dense in $\Diff^1_P ([0,1])$.

To see that $\mathcal{S}_n$ is open, fix $f \in \mathcal{S}_n$, and let $a,b$ be points in $(0,1/n)$ for which 
$f(x)-x$ is positive (resp. negative) at $a$ (resp. $b$). If $g$ is a $C^1$-small perturbation of $f$, 
then $g(x)-x$ is still positive (resp. negative) at $a$ (resp. $b$). 
This shows that $\mathcal{S}_n$ contains a neighborhood of $f$. 
Since $f \in \mathcal{S}_n$ was arbitrary, this implies that $\mathcal{S}_n$ is open.

To see that $\mathcal{S}_n$ is dense, recall from Lemma \ref{l:G1} that the subset $G_1$ of $\Diff^1_P([0,1])$ made of  the diffeomorphisms 
with finitely many fixed points is dense. We claim that each $f \in G_1$ can be approximated in $C^1$ topology by elements in $\mathcal{S}_n$. 
To do this, fix $\varepsilon > 0$, and for each $\delta$ consider the affine map $\varphi_{\delta}$ that sends $[0,1]$ to $[\delta,1]$. Fix any $C^1$ 
diffeomorphism $h_{\delta}$ of $[0,\delta]$ into itself such that:
\begin{itemize}
\item it is parabolic at the endpoints;
\item it has exactly one interior fixed point, which is also parabolic and {\em topologically transversal}  
(i.e. of type 2, in the terminology of the proof of Lemma \ref{l:G2});
\item it is $\varepsilon$-close to the identity in $C^1$ topology.
\end{itemize}
Choosing appropriately $\delta < 1/n$, the diffeomorphism $f_{\delta}$ defined as
$$f_{\delta} (x)
=
\begin{cases} 
          \varphi_{\delta} f \varphi_{\delta}^{-1} (x) & \mbox{if } x \in [\delta,1] , \\
           h_{\delta} (x) & \mbox{if } x \in  [0,\delta], \\
\end{cases}
$$
is $\varepsilon$-close to $f$ in $C^1$ norm. Since it has a topologically transversal fixed point in $(0,1/n)$, 
this element $f_{\delta}$ belongs to $\mathcal{S}_n$. This concludes the proof.
\end{proof}

\begin{proof}[First proof of Theorem \ref{t:generic}]
By the preceding lemma, for a generic $f \in \Diff^1_P([0,1])$, the function $x \mapsto f(x)-x$ changes its sign 
infinitely many times in every neighborhood of the origin. 
We claim that such an $f$ cannot be $C^1$-conjugated to a nontrivial power of itself. Indeed, if $hfh^{-1}=f^k$ for a certain 
$h \in \Diff_+^1([0,1])$ and $k\geq 2$ then, by Theorem \ref{t:CC}, the conjugating diffeomorphism $h$ fixes all the infinitely 
many components of $f$. By Theorem \ref{t:BMNR}, inside each of these components there is a point at which its  
derivative equals $k$. Since there is a sequence of components converging to the origin, this implies that $Dh (0) = k$. 
However, since $h$ also has a sequence of fixed points converging to the origin, one necessarily has $Dh (0) = 1$. 
This contradiction concludes the proof.
\end{proof}


\section{A second proof of Theorem \ref{t:generic}}
\label{s:generic}

Our second proof of Theorem \ref{t:generic} gives more information than the first one. It relies on an adaptation of tools introduced 
in \cite{BCW09} (and a preliminary version of it graciously communicated to us by one of the authors) to prove that the $C^1$ generic 
diffeomorphism $f$ of any compact manifold has a trivial centralizer (i.e. it is reduced to the cyclic group generated by $f$). 
Observe that such a diffeomorphism cannot be conjugated to a nontrivial power $f^k$. Indeed, if $ h f h^{-1} = f^k$,  
then $(h^{-1} f h)^k = f$, hence $f$ would have a $k^{th}$-root in its centralizer (namely, $h^{-1}fh$), 
which would imply that the centralizer of $f$ is not generated by $f$.

However, having a trivial centralizer is not a generic nor even dense property in $\Diff_P^1([0,1])$. Indeed, as soon as $f\in \Diff_P^1([0,1])$ 
has an interior (parabolic) fixed point $p$ and differs from the identity on both sides of $p$, the map $\hat{f}$ defined as $f$ on $[0,p]$ and $f^2$ on $[p,1]$, 
for example, is $C^1$ (thanks to the parabolicity of $p$), obviously commutes with $f$, and is not a power of $f$. Now take such an $f$ with the additional 
property that $f(x)-x$ is positive on $(0,p)$ and negative on $(p,1)$, say. Every $C^1$-small perturbation $g \in \Diff_P^1([0,1])$ of $f$ will satisfy that 
$g(x)-x$ takes different signs, hence it has an interior fixed point. Consequently, $g$ has a nontrivial centralizer. By this procedure we 
have thus detected an open set of diffeomorphisms whose centralizers are nontrivial.

\medskip

The above discussion motivates the following definition: 

\begin{defn}
We say that $f \in \Diff_P^1([0,1])$ is {\em locally trivially centralized} if $f$ is not the identity and, for every component $I$ of $f$, the 
centralizer $\mathrm{Z}_{I,f}^1$ in $\Diff^1_+(I)$ of the diffeomorphism of $I$ induced by $f$ is trivial. 
\end{defn}

\medskip

We denote by $\LTC$ the subset of $\Diff_P^1([0,1])$ made of the elements that are locally trivially centralized. 
The claim in the first paragraph of this section has a natural extension to this context.

\begin{lem} 
No element in $\LTC$ is $C^1$-conjugated to a nontrivial power of itself.
\end{lem}

\begin{proof} 
Let $f \in \Diff_P^1 ([0,1])$ be nontrivial and let $I$ one of its components. If $f$ is conjugate to a nontrivial 
power of itself, say $ h f h^{-1} = f^k$,  then $(h^{-1} f h)^k = f$, thus $f$ has a $k^{th}$-root, namely, $g := h^{-1}fh$. 
It is not hard to see that $I$ is also a component for $g$. Indeed, $g$ cannot have a fixed point inside (otherwise $f$ would 
also have one), and has to fix the endpoints of $I$ (otherwise $f$ would not fix them). Therefore, the diffeomorphism of $I$ 
induced by $g$ is a $k^{th}$-root of the one induced by $f$, which shows that $f$ is not locally trivially centralized.
\end{proof}

Due to the preceding lemma, Theorem~\ref{t:generic} directly follows from the next more general result.

\begin{thm}
\label{t:almost-trivial}
The set $\LTC$ is residual in $\Diff^1_P([0,1])$.
\end{thm}

As we next show, the tools introduced in \S \ref{s:distorted} allow to easily show that $\LTC$ is dense in $\Diff^1_P([0,1])$. 
However, proving that it is a residual set is much harder. The key  idea will be explained in \S \ref{s:idea}. This idea, as 
well as the subsequent proof of Theorem \ref{t:almost-trivial}, are extremely close to the preliminary version of \cite{BCW09} 
already alluded to. The main difference is that, as in \S~\ref{s:distorted}, hyperbolic fixed points are not allowed in our context.

\medskip

\begin{prop}
\label{p:almost-trivial}
The set $\LTC$ is dense in $\Diff^1_P([0,1])$.
\end{prop}

\begin{proof} By Lemma \ref{l:G3}, the subset of $\Diff^1_P([0,1])$ made of $C^\infty$-flowable diffeomorphisms with a finite 
number of fixed points is dense in $\Diff^1_P([0,1])$. It thus suffices to prove that any such diffeomorphism $f$ can be approached by an element 
of $\LTC$. Assume first that $f$ has no interior fixed points. Now recall that $C^1$-flowability of such a diffeomorphism, which corresponds to the triviality 
of its Mather invariant (cf. \S \ref{s:back}), is ``broken'' by \emph{any} (nontrivial) perturbation supported in a single fundamental domain. The 
resulting diffeomorphism $g$ then has an infinite cyclic $C^1$-centralizer. Moreover, one can actually prescribe the (nontrivial) Mather invariant 
$M_g$ of $g$ among the (class up to pre and post composition by translations of) diffeomorphisms $\tilde \phi$ of $\R$ commuting with the 
unit translation and having intervals of fixed points (all of this directly follows from Lemma \ref{l:modify}). In particular, one can ensure that 
$M_g$ does not commute with the translation by $1/k$ for any $k \geq 2$. Again by \S \ref{s:back}, the latter implies that $g$ has no $k^{th}$ 
root in $\Diff^1_+([0,1])$ for any $k \geq 2$, and thus its $C^1$-centralizer is trivial. Since $g$ has no interior fixed point, this implies that it 
belongs to $\LTC$. This concludes the proof in the case where $f$ has no interior fixed point.

Now if $f$ has (a finite number of) interior fixed points (as at the beginning of the proof), one just needs to 
perform the above procedure on every component of $f$, and the resulting diffeomorphism remains 
$C^\infty$ because it is smooth inside each component and nothing has been changed near the fixed points. 
\end{proof}


\subsection{Genericity of $\LTC$ in $\Diff^1_P([0,1])$: a key idea for the proof}
\label{s:idea}

The key idea behind the proof of Theorem \ref{t:almost-trivial} comes from the fact that large centralizers allow controlling dynamical distortion, that 
is, obtaining  a uniform bound for the quotient of the derivatives of iterates along different orbits. This is the content of Lemma \ref{l:stupid} below. 
However, as we will see in Proposition \ref{p:residual}, absence of control of dynamical distortion is generic among $C^1$ diffeomorphisms 
with no hyperbolic fixed points.

\begin{lem}
\label{l:stupid}
Let $f\in\Diff^1_+([0,1])$ be a nontrivial element and $I$ a component of it. For every $x \in I^{\mathrm{o}}$ and all $g\in\mathrm{Z}_{I,f}^1$, one has 
$$\sup_{n\in\N}|\log Df^n(x)-\log Df^n(g(x))|<+\infty.$$
\end{lem}

\begin{proof} On the interval $I$ we have $gf^n=f^ng$, which by the chain rule gives 
$$\log Dg(f^n(x))+\log Df^n(x) = \log Df^n (g(x))+\log Dg(x),$$
thus 
$$\log Df^n(x) - \log Df^n (g(x)) = \log Dg(x) - \log Dg(f^n(x)).$$
As a consequence, for all $n \in \N$, 
$$|\log Df^n(x)-\log Df^n(g(x))|$$ 
is bounded from above by $\, 2 \| \! \log Dg \|_{\infty}$.
\end{proof}

Before stating the genericity of absence of control of dynamical distortion, we introduce the notation we will use concerning what is 
called {\em the unbounded distortion property} in \cite{BCW09}.

\begin{defn} We denote by $\UD$ the set of nontrivial elements $f \in \Diff_P^1([0,1])$ satisfying the following property: For all $x,y$ in the same 
component of $[0,1]\setminus\Fix(f)$ and such that $O_f(x)\cap O_f(y)=\varnothing$, one has 
$$\sup_{n\in\N}|\log Df^n(x)-\log Df^n(y)|=+\infty.$$
\end{defn}

\begin{prop}
\label{p:residual}
The set $\, \UD$ is residual in $\Diff_P^1([0,1])$.
\end{prop}
The proof of this proposition will be the content of \S \ref{ss:Gdelta} and \ref{ss:density}.
Now we show how Theorem~\ref{t:almost-trivial} (hence Theorem \ref{t:generic}) follows from it.

\begin{proof}[Proof of Theorem \ref{t:almost-trivial} from Proposition \ref{p:residual}] 
It suffices to show that every $f\in \UD$ is locally trivially centralized. Let $I$ be a component of $f$, let $x \in I^{\mathrm{o}}$ 
and let $g \in \mathrm{Z}^1_{I,f}$. Assume that the point $z:=g(x)$, 
which belongs to $I^{\mathrm{o}}$, does not belong to $O_f(x)$, and hence $O_f(x)\cap O_f(z)=\varnothing$. 
Then, by the definition of $\UD$, we have 
$$\sup_{n\in\N}|\log Df^n(x)-\log Df^n(z)|=+\infty.$$
However, this contradicts Lemma \ref{l:stupid}.

We conclude that  $g(x)$ belongs to $O_f(x)$, thus there exists $n(x)\in\Z$ such that $g(x)=f^{n(x)}(x)$. 
Since $g$ is continuous and $O_f(y)$ is discrete for every $y \in I^{\mathrm{o}}$, the map $x \mapsto n(x)$ 
is easily seen to be constant on $I$, say equal to $n$. This implies that $g = f^n$ (on $I$), which concludes the proof.
\end{proof}

Now, to prove Proposition \ref{p:residual}, given $f$ and \emph{some} pair $x,y$ in the interior of the same component $I$ of 
$f$  such that $O_f(x)\cap O_f(y)=\varnothing$, it is rather elementary to find a small perturbation $g$ of $f$ that satisfies 
$$\sup_{n\in\N}|\log Dg^n(x)-\log Dg^n(y)|=+\infty.$$ 
However, it is more difficult to find a single perturbation that satisfies this distortion phenomenon for \emph{each} such pair $x,y$ (this is 
equivalent to proving the density of $\UD$ in $\Diff_P^1([0,1])$). Finally, proving that $\UD$ is residual is even more complicated. These two 
issues will be the content of \S \ref{ss:Gdelta} and \S \ref{ss:density} below.


\subsection{Proof of Proposition \ref{p:residual}, part I: definition of a dense $G_\delta$-set} 

\label{ss:Gdelta}

Let us denote by $\cI$ be the countable family of closed subintervals of $(0,1)$ with rational endpoints. 
Similarly, for every $\Delta\in \cI$, let $\cI_\Delta$ denote the countable set of closed subintervals of $\Delta$ having rational endpoints. 
Given $\Delta \in\cI$, we consider the open subset of $\Diff_P^1 ([0,1])$ defined by 
$$O(\Delta):=\{f\in \Diff_P^1([0,1]), f(\Delta)\cap \Delta=\varnothing\}.$$
Given two disjoint compact subintervals $J_1,J_2$ of $\Delta$ and 
$M>0$,
let 
$$O(\Delta, J_1,J_2,M):=\left\{f\in O(\Delta)\; | \; \exists n\in\N:\forall x_1\in J_1, \forall x_2\in J_2, |\log Df^n(x_2)-\log Df^n(x_1)|>M\right\}.$$
We immediately state a crucial lemma for the proof of Proposition \ref{p:residual}, the proof of which will be given later. 
\begin{lem}
\label{l:dense} 
For every $\Delta, J_1, J_2, M$ as above, the set 
$$O(\Delta, J_1,J_2,M)\cup \Int(\Diff_P^1([0,1])\setminus O(\Delta))$$ 
is open and dense in $\Diff^1_P([0,1])$.
\end{lem}

\begin{proof}[Proof of Proposition \ref{p:residual} assuming Lemma \ref{l:dense}] 
Our assumption  implies that the (countable) intersection 
$$\cR := \bigcap (O(\Delta, J_1,J_2,M)\cup \Int(\Diff_P^1([0,1])\setminus O(\Delta)))$$
taken over all $\Delta\in \cI$, all $ J_1, J_2$ in $\cI_{\Delta}$ with $J_1\cap J_2 = \varnothing$, and all $M\in\N$, 
is a dense $G_\delta$-set in $\Diff_P^1([0,1])$. We next show that $\cR$ is contained in $\UD$, which implies that the latter set is residual.

Let $f\in\cR$ and let $x, y$ in the interior of a component $I$ of $f$ be such that $O_f(x)\cap O_f(y)=\varnothing$. 
Let $L$ be a fundamental domain of $f$ containing $x$ in its interior, and let $m\in \Z$ be such that $f^m(y) \in L$. 
By slightly moving $L$, we can assume that $f^m(y)$ is also in its interior. Let $\Delta$ be a compact subinterval of $L$ still 
containing $x$ and $f^m (y)$ in its interior. By construction, $f$ lies in $O(\Delta)$, and thus $f\notin \Int(\Diff_P^1([0,1])\setminus O(\Delta))$.

Now choose $J_1,J_2\in \cI_{\Delta^{\mathrm{o}}}$ such that $J_1\cap J_2 = \varnothing$ and $x \in J_1$, $y \in J_2$. 
Since $f$ lies in $\cR$ and $f\notin \Int(\Diff_P^1([0,1])\setminus O(\Delta))$, it must belong to $O(\Delta, J_1, J_2,M)$ for every $M > 0$. 
Thus, $f$ belongs to the intersection of these sets over all $M\in\N$.   
In particular, for each $M > 0$ there exists $n\in\N$ such that 
$$|\log Df^n(x)-\log Df^n(f^m(y))|>M.$$
As a consequence, 
$$\sup_n|\log Df^n(x)-\log Df^n(f^m(y))|=+\infty.$$
Since, by the chain rule,  
$$|\log Df^n(f^m(y)) - \log Df^n (y)| = | \log Df^m (f^n (y))- \log Df^m (y) | \leq 2m \| \log Df\|_{\infty},$$
this implies that 
$$\sup_n|\log Df^n(x)-\log Df^n(y)| = +\infty.$$
Therefore, $f$ belongs to $\UD$, as claimed.
\end{proof}

We are thus left to proving Lemma \ref{l:dense}. This will be a consequence of the next crucial lemma, 
which will be proved in \S \ref{ss:density}.

\medskip 

\begin{lem}
\label{l:ddense}
For every $\Delta, J_1, J_2, M$ as in Lemma \ref{l:dense}, the set $O(\Delta, J_1,J_2,M)$ is open and dense in $O(\Delta)$.
\end{lem}

\medskip

\noindent{\em Proof of Lemma \ref{l:dense} assuming Lemma \ref{l:ddense}.}
That $O(\Delta, J_1,J_2,M)$ is open trivially follows from its definition. 
This obviously implies that the set 
$$O(\Delta, J_1,J_2,M)\cup \Int(\Diff_P^1([0,1])\setminus O(\Delta))$$ 
is also open. To show that this set is also dense in 
$\Diff^1_P([0,1])$, fix an element $f$ therein. There are three possibilities: 
\begin{itemize}
\item $f(\Delta)\cap\Delta$ has nonempty interior: then this still holds for any small enough perturbation of $f$, 
hence $f$ belongs to $\, \Int(\Diff_P^1([0,1])\setminus O(\Delta))$;
\item $f(\Delta)\cap\Delta$ is empty: this means that $f$ belongs to $O(\Delta)$ and hence, by Lemma \ref{l:ddense}, 
it can be approached by elements of $O(\Delta, J_1,J_2,M)$;
\item $f(\Delta)\cap\Delta$ is nonempty but has no interior: this means that $f(\Delta)\cap\Delta$ is made of a single point and then, 
by a well-chosen arbitrarily small perturbation $g$, one can ensure that $g(\Delta)\cap\Delta$ has some interior. As a consequence, 
$f$ can be approached by elements of $\Int(\Diff_P^1([0,1])\setminus O(\Delta))$.
 $\hfill\square$
\end{itemize}


\subsection{Proof of Proposition \ref{p:residual}, part II: density of unbounded dynamical distortion}
\label{ss:density}

Here, we prove the key Lemma \ref{l:ddense}.  
Let $\Delta\in \cI$, let $J_1, J_2$ be compact disjoint subintervals of $\Delta$, and let $M>0$. Recall that
$$O(\Delta, J_1,J_2,M):=\left\{f\in O(\Delta)\; | \; \exists n\in\N,\forall x_1\in J_1, x_2\in J_2, |\log Df^n(x_2)-\log Df^n(x_1)|>M\right\}.$$
From this definition, it immediately follows that $O(\Delta, J_1,J_2,M)$ is open. We next focus on its density in $O(\Delta)$. We let 
$f\in O(\Delta)$ and $\eps>0$. We want to construct $g\in O(\Delta, J_1,J_2,M)$ such that $\| \log Df - \log Dg \|_{\infty} <\eps$. 

\paragraph{Preparation.} Since the interval $\Delta$ satisfies $f(\Delta)\cap \Delta=\varnothing$, it must lie in some fundamental domain 
$L$ of $f$, itself lying in some component $I = [a,b]$ of $f$. To fix ideas, we will assume that $f$ ``pushes everything to the right'' on $I$.

Mimicking the proof of Lemma \ref{l:G2}, one can make $f$ of class $C^2$ on a left neighborhood of $b$ by an arbitrarily $C^1$-small perturbation 
supported in a left neighborfood of $b$ away from $\Delta\cup f(\Delta)$, all of this without adding new fixed points. We can then fix a left neighborhood 
$U$ of $b$, again away from $\Delta\cup f(\Delta)$, such that $V(f;U) < \eps / 4$. 
(and in particular $|\log Df|<\eps/4$ on $U,$ since $b$ is parabolic). We will still denote the resulting perturbed diffeomorphism by~$f$.

Fix a sufficiently large  $m\in\N$ such that $f^m(\Delta)\subset U$. Of course, $f^m(J_1)$ and $f^m(J_2)$ are still disjoint (and lie in $U$). 
Choose points $p_1< p_2 <p_3<p_4<p_5$ in $f^m(\Delta)$ such that $f^m(J_1)\subset (p_1,p_3)$ and $f^m(J_2)\subset (p_3,p_4)$.

Consider a diffeomorphism $h$ of the interval $f^m(\Delta)$ such that:
\begin{itemize}
\item $h=\id$ outside $(p_1,p_5)$;
\item the set of fixed points of $h$ in $(p_1,p_5)$ coincides with $\{p_2, p_3, p_4\}$;
\item there exist $\mu_1, \mu_2$ such that $0<Dh(p_2)<\mu_1<\mu_2<Dh(p_4)<1$;
\item $\| \log D h  \|_{\infty} < \eps / 4$.
\end{itemize}

\paragraph{Construction of a sequence $f_n\in O(\Delta)$ tending to $f$.} For every $n>0$, we define $f_n$ as follows: 
$$f_n
=
\begin{cases} 
      f \circ (f^j h f^{-j}) & \mbox{on } f^{j+m}(\Delta), \, j \in \{0,1,\ldots,n\} , \\
      f & \mbox{elsewhere}. \\
\end{cases}
$$

It is immediate that the $f_n$ are $C^1$ diffeomorphisms of $[0,1]$ that still displace $\Delta$ from itself (since we did not change $f$ 
in a neighborhood of $\Delta\cup f(\Delta)$). To prove that $f_n$ belongs to $O(\Delta)$, we must check that $f_n$ belongs to 
$\Diff^1_P([0,1])$, which amounts to showing that it has only parabolic fixed points. This follows from the next stronger claim:

\medskip

\begin{claim}
\label{c:claim2}
The map $f_n$ has the same fixed points as $f$ and coincides with $f$ near its fixed points.
\end{claim}

\begin{proof}
The only place where $f_n$ differs from $f$ are the intervals $f^{j+m}(\Delta)$. Each of these is fixed by $f^j h f^{-j}$ and sent 
into a disjoint interval by $f$, hence by $f_n$, and thus there is no new fixed point of $f_n$ therein.
\end{proof}

We next prove two more claims that, putted together, allow to conclude the proof of the density of $O(\Delta, J_1,J_2,M)$ 
in $O(\Delta)$.

\medskip

\begin{claim}
\label{c:claim3}
For every $n$, one has $\| \log D f - \log D f_n  \|_{\infty} < \eps$.
\end{claim}

\begin{proof} 
Outside the intervals $f^{j+m}(\Delta)$, with $j \in \{0,1,\ldots,n\}$, the diffeomorphism $f_n$ coincides with $f$. 
On the intervals $f^{j+m}(\Delta)$, one has $f_n = f \circ (f^j h f^{-j})$, hence 
\begin{align*}
\log Df - \log Df_n 
&= \log Df - \log D(f\circ f^jhf^{-j})\\
&=  \log Df - \log Df(f^jhf^{-j})) - \log D(f^jhf^{-j}).
\end{align*}
Since $f^{j+m}(\Delta)\subset U$, this implies that, still on  $f^{j+m}(\Delta)$, 
$$|\log Df - \log Df_n|\le |\log Df|+ |\log Df(f^jhf^{-j}))| + |\log D(f^jhf^{-j})|<\frac\eps4+\frac\eps4+|\log D(f^jhf^{-j})|.$$
Now, for $x\in f^{j+m}(\Delta)$,
\begin{align*}
|\log D(f^jhf^{-j})(x)|
&=|\log Df^j (hf^{-j}(x))-\log Df^j(f^{-j}(x))+\log Dh(f^{-j}(x))|\\
&\le \frac\eps4+ \sum_{k=0}^{j-1} |\log Df(f^khf^{-j}(x))-\log Df(f^kf^{-j}(x))|\\
&\le \frac\eps4 + V (f;U)\\
&< \frac{\eps}{2},
\end{align*}
where the second inequality comes from the fact that $hf^{-j}(x)$ and $f^{-j}(x)$ belong to the same 
fundamental domain of $f$ and that this, as well as all of its forward images by $f$, lies in $U$. 
In the end, we get \, $\| \! \log D f - \log D f_n \|_1 < \eps$, as announced.
\end{proof}

\medskip

\begin{claim}
\label{c:claim4}
For $n$ large enough, $f_n$ belongs to $O(\Delta,J_1,J_2,M)$.
\end{claim}

\begin{proof}
More precisely, we claim that for a large enough $n$, for all $x_1\in J_1$ and $x_2\in J_2$  one has 
$$|\log Df_n^{n+m}(x_2)-\log Df_n^{n+m}(x_1)| > M.$$ 
To prove this, we will show that 
$$\frac{D(f_n)^{n+m}(x_2)}{D(f_n)^{n+m}(x_1)} > e^M.$$

Note that for every $q_1\in f^m(J_1)\subset (p_1,p_3)$ and $q_2\in f^m(J_2)\subset (p_3,p_4)$, 
the second and third properties of $h$ give, for $n$ larger than some fixed $N_0$, 
\begin{equation}
Dh^{n+1}(q_1)<\mu_1^n \quad \text{ and } \quad  Dh^{n+1}(q_2)>\mu_2^n.
\end{equation}
Let us now fix $\, C\ge 1 \,$ such that 
\begin{equation}
\label{e:C} 
\frac{Df^m(y)}{Df^m(x)}\le C \quad \mbox{ for all }  x,y \mbox{ in } \Delta 
\end{equation}
and $\, e^{V(f; U)} \le C. \,$ Now choose $n\ge N_0$ such that 
\begin{equation}
e^M < \frac1{C^2}\left(\frac{\mu_1}{\mu_2}\right)^n.
\end{equation}
By definition, for $q\in f^m(\Delta)$, we have 
$$(f_n)^n(q) = \left((f^{n+1}h f^{-n})\circ (f^nhf^{-n+1})\circ \cdots\circ fh\right)(q) =(f^{n+1}\circ h^{n+1})(q).$$
Given $x_1\in J_1$ and $x_2\in J_2$, let $q_1:=f^m(x_1)=(f_n)^m(x_1)$ and $q_2:=f^m(x_2)=(f_n)^m(x_2)$. 
Then, according to \eqref{e:C}, 
\begin{align*}
\frac{D(f_n)^{m+n}(x_2)}{D(f_n)^{m+n}(x_1)}
&= \frac{D(f_n)^n((f_n)^m(x_2))D(f_n)^m(x_2)}{D(f_n)^n((f_n)^m(x_1))D(f_n)^m(x_1)}\\
&= \frac{D(f_n)^n(q_2)Df^m(x_2)}{D(f_n)^n(q_1)Df^m(x_1)}\\
&\ge C^{-1} \frac{D(f_n)^n(q_2)}{D(f_n)^n(q_1)} \\
&= C^{-1} \frac{Df^{n+1}(h^{n+1}(q_2))Dh^{n+1}(q_2)}{Df^{n+1}(h^{n+1}(q_1))Dh^{n+1}(q_1)}
\end{align*}
Now letting $p_1 := h^{n+1}(q_1)\in U$ and $p_2 := h^{n+1}(q_2)\in U$, we obtain 
\begin{align*}
\left|\log\frac{Df^{n+1}(p_2)}{Df^{n+1}(p_1)}\right|
&=\left|\log Df^{n+1}(p_2)- \log Df^{n+1}(p_1)\right|\leq V(f;U),
\end{align*}
and therefore
$$\frac{Df^{n+1}(p_2)}{Df^{n+1}(p_1)}\geq e^{-V (f;U)}\geq C^{-1}.$$
We thus finally obtain 
$$\frac{D(f_n)^{m+n}(x_2)}{D(f_n)^{m+n}(x_1)}\ge C^{-2} \left(\frac{\mu_2}{\mu_1}\right)^n>e^M,$$
which concludes the proof.
\end{proof}


\section{A proof of Theorem \ref{t:higher}}
\label{s:higher}

In order to prove Theorem \ref{t:higher}, we first introduce a powerful tool, the \emph{asymptotic variation}, which both detects undistortion 
in regularity higher than $C^{1+bv}$ and non-$C^1$-conjugacy to a power. Although defined independently, it is closely related to the 
Mather invariant. 

\medskip

Recall that for $f \in \Diff^{1+bv}_+([0,1])$ we denote $V(f)$ the total variation of the logarithm of its derivative. 
The chain rule $D(fg) = Dg \cdot Df(g)$ easily gives the general subadditive inequality 
$$V(fg) \leq V(f) + V(g)$$ 
for every $f,g$ in $\Diff^{1+bv}_+([0,1])$. This implies in particular 
that the sequence $V(f^n)$ is subadditive. By the classical Fekete's lemma, the sequence 
$V(f^n) / n$ converges. 
We denote the limit by $V_{\infty}(f)$ and we call it the {\em asymptotic variation}:
$$V_{\infty} (f) := \lim_{n \to \infty} \frac{V(f^n)}{n}.$$
This has many interesting properties, among which the following will be crucial for our purposes: 
\begin{enumerate}
\item For all $k \geq 1$ and all $f \in \Diff^{1+bv}_+([0,1])$, one has $V_{\infty}(f^k) = k \cdot V_{\infty}(f)$; 
\item If $f$ is a distortion element of any group $G$ contained in $\Diff^{1+bv}_+([0,1])$, then $V_{\infty} (f) = 0$; 
\item A diffeomorphism $f \in \Diff_+^{2} ([0,1])$ has vanishing asymptotic variation if and only if all its fixed points are parabolic and it is $C^1$-flowable;
\item If $f_n$ is a sequence of $C^{1+bv}$ diffeomorphisms of $[0,1]$ converging to $f$ in the $C^{1+bv}$ sense  
and neither $f$ nor any $f_n$ has interior fixed points, then $V_{\infty} (f_n)$ converges to $V_{\infty} (f)$;
\item If $f \in \Diff_P^{1+bv} ([0,1])$ has no interior fixed point and $h \in \mathrm{Diff}_+^1 ([0,1])$ is such that the 
conjugate $g:= hfh^{-1}$ is a $C^{1+bv}$ diffeomorphism, then $V_{\infty} (g) = V_{\infty} (f)$.
\end{enumerate}

The (homogeneity) Property 1 above follows from the very definition of $V_{\infty}$, and Property 2 follows from the subadditive 
property of $V$ and the definition of distortion element.  For Properties 3, 4 and 5, we refer to \cite{EN21} 
(cf. Theorems A, C and D therein, respectively). Note that, given 
$f\in\Diff^{1+bv}([0,1])$ and a component $I$ of $f$, one can similarly define the asymptotic variation of $f$ restricted to 
$I$, that we denote by $V_\infty(f;I)$. With this notation, we have the following localization property, which 
corresponds to \cite[Lemma 1.1]{EN21}: 
\begin{equation}\label{e:localization}
V_{\infty} (f) = \sum_{I \in \mathcal{I}} V_{\infty} (f;I).
\end{equation}
Here $\mathcal{I}$ denotes the family of components of $f$. 
To all these properties, we add the following proposition, the proof of which will be given in \S \ref{s:Vconj}.

\medskip

\begin{prop}
\label{p:Vconj}
If $f \in \Diff_P^{2} ([0,1])$ has nonvanishing asymptotic distortion, then it cannot be $C^1$-conjugated to a nontrivial power.
\end{prop}

Let us now fix $\ell \geq 2$ (which may be an integer number or infinite)  
and consider the set $\cV_{\ell}:=\{f\in\Diff^{\ell}_P([0,1]): V_\infty(f)>0\}$. From Property~2 and Proposition \ref{p:Vconj}, 
it immediately follows that $\cV_{\ell}$ lies inside the subset  $\mathcal{T}_{\ell}$ of $\mathrm{Diff}^{\ell}_p ([0,1])$ made 
of the elements which are neither $C^{\ell}$-distorted nor $C^1$ conjugated to some nontrivial power of themselves. To prove 
Theorem \ref{t:higher}, it thus suffices to show the next two statements, which will be treated in subsequent subsections.

\medskip

\begin{prop}
\label{p:notopen}
The set $\mathcal{T}_{\ell}$ is not open in $\mathrm{Diff}^{\ell}_p ([0,1])$.
\end{prop}

\begin{prop}
\label{p:Vdense}
The interior 
of $\cV_{\ell}$ is dense in $\mathrm{Diff}^{\ell}_p ([0,1])$.
\end{prop}


\subsection{Two proofs of Proposition \ref{p:Vconj}}
\label{s:Vconj}

There are at least two ways of proving Proposition \ref{p:Vconj}. To begin with, 
note that, if $V_{\infty}(f) > 0$, then the localization property (\ref{e:localization}) 
gives us a component $I$ of $f$ for which $V_\infty (f;I)>0$.

Now assume that $f$ is conjugated by some $C^1$ diffeomorphism $h$ of $[0,1]$ to a power $f^k$, with $k\ge 2$. 
Then Properties 1 and 5 above (which easily extend to diffeomorphisms defined in different intervals) give 
$$k \cdot V_{\infty} (f;I) = V_{\infty} (f^k; I) = V_{\infty} (hfh^{-1};I) = V_{\infty} (f; I),$$
which is a contradiction.

For a second proof, note that, on the first hand, the relation $hfh^{-1} = f^k$ gives $(h^{-1} f h)^k = f$, 
which means that $f$ has a nontrivial root.  
On the other hand, by Property~3 above, the fact that $V_{\infty} (f;I) > 0$ implies that it is not $C^1$-flowable on $I$, which means 
that its Mather invariant (to the diffeomorphism of $I$ induced by $f$) is nontrivial. However, by Remark \ref{r:mather-root}, this is 
incompatible with $f$ having a nontrivial root.


\subsection{Proof of Proposition \ref{p:notopen}}

This follows almost directly from a construction in \cite[\S 6]{EN21} of a sequence $(f_n)$ of $C^\infty$-flowable diffeomorphisms 
without interior fixed points, parabolic at the boundary  (thus with vanishing asymptotic distortion, according to Property 3 above), 
which converge to a diffeomorphism $f$ with a single interior fixed point and positive asymptotic variation on both sides. Indeed, 
this construction is very flexible. In particular, one can impose the germs of the $f_n$'s at $0$ and $1$ to be $C^\infty$-conjugated 
to their square, which, together with the flowability property, implies that $f_n$ is $C^\infty$-conjugated to its square for every $n$ 
(cf. Proposition \ref{l:criterion}), and in particular $C^\infty$-distorted. The limit $f$, however, is not $C^{\ell}$-distorted in 
virtue of Property 2. Moreover, the nonvanishing of its asymptotic distortion prevents it to be $C^1$-conjugated to a 
nontrivial power of itself, as it was argued above.


\subsection{Proof of Proposition \ref{p:Vdense}}
\label{s:Vdense}

Before passing to the proof, it is worth mentioning that we deal with the interior of 
$\cV_{\ell}$ rather than $\cV_{\ell}$ itself since, as it was recalled above, the latter set is not open. 

Let us denote by $\mathcal{U}_{\ell}$ the set of diffeomorphisms $f$ of $\mathrm{Diff}^{\ell}_p ([0,1])$ 
which have a component $[a,b]$ such that $V_{\infty} (f;[a,b]) > 0$ and, if $a$ (resp. $b$) lies in the interior of $[0,1]$,  
then $a$ (resp. $b$) is a topologically transversal fixed point of $f$. Proposition \ref{p:Vdense} then follows from:

\medskip

\begin{prop}
\label{p:D1}
Every element of $\mathcal{U}_{\ell}$ lies in the interior of $\cV_{\ell}$.
\end{prop}

\begin{prop}
\label{p:D2}
The set $\mathcal{U}_{\ell}$ is dense in $\mathrm{Diff}^{\ell}_p ([0,1])$.
\end{prop}

\medskip

\begin{proof}[Proof of Proposition \ref{p:D1}] 
Fix $f\in \mathcal{U}_{\ell}$, and let $I = [a,b]$ be a component of $f$ satisfying the properties of the definition of $\mathcal{U}_{\ell}$. 
If $(f_n)$ is a sequence in $\mathrm{Diff}^{\ell}_p ([0,1])$ that converges to $f$ then, by transversality, each $f_n$ has a component $[a_n,b_n]$  
close to $[a,b]$ for large enough $n$. By the (continuity) Property 4 above (which easily extends to diffeomorphisms defined on different 
though very close intervals), one has that $V_{\infty} (f_n; [a_n,b_n])$ converges to $V_{\infty} (f;[a,b])$. 
We thus conclude that $V_{\infty} (f_n; [a_n, b_n])$, and hence $V_{\infty} (f_n)$, is positive for large enough $n$. 
This shows that $f$ lies in the interior of $\cV_{\ell}$. 
\end{proof}

\medskip

\begin{proof}[Proof of Proposition \ref{p:D2}] 
We first show that the subset $G_{1,\ell}$ of $\Diff_P^{\ell}([0,1])$ formed by the diffeomorphisms with finitely many fixed 
points is dense. The argument is very close to that of Lemma~\ref{l:G1}, with a slight change. 

First assume that $\ell$ is finite. Given $f\in \Diff_P^{\ell} ([0,1])$ 
and $\eps>0$, consider the family $\mathcal{I}$ formed by the maximal intervals $I$ whose boundary points are fixed points 
of $f$ and where $\| (f-\id)|_I \|_{\ell} \leq \eps / 3$. On each interval $I$ of $\mathcal{I}$, replace $f$ by any $C^{\ell}$ diffeomorphism 
that has the same Taylor series expansion as $f$ up to order $\ell$ at the endpoints, has at most one fixed point in the interior and is $\eps/2$-close 
to the identity in $C^{\ell}$ norm. (Note that one fixed point inside $I$ is necessary if, at both endpoints of $I$, the map $f$ is simultaneously 
topologically contracting or topologically expanding.) 
Let $\hat{f}$ be the new $C^{\ell}$ diffeomorphism thus obtained. By construction, $\Fix(\hat f)\subset \Fix(f)$. 
Moreover, $\hat{f}$ coincides with $f$ outside the intervals in $\mathcal{I}$, whereas on each $I \in \mathcal{I}$ we have
$$\| (\hat {f} - f)|_I\|_{\ell} \leq \| (\hat{f} - \id)|_I \|_{\ell} + \| ( \id - f )|_{\ell}  \|_1 \leq \frac{\varepsilon}{3} + \frac{\varepsilon}{2} < \varepsilon.$$
This shows that $\hat{f}$ is $\varepsilon$-close to $f$ in $C^{\ell}$ norm. 
Finally, we claim that $\Fix(\hat f)$ is finite. Indeed, assume by contradiction that this is not the case and let $p$ be an accumulation 
point of $\Fix(\hat f)$, say from the right. Then there exists an interval $[p,q]$ containing infinitely many fixed points of $\hat{f}$, hence 
of $f$. This implies that $f$ is $\ell$-tangent to the identity at $p$. Thus, by slightly approaching $q$ to $p$ if necessary, we may 
assume that $\| (f-\id)|_{[p,q]} \|_{\ell} \leq \varepsilon / 3$. However, this is in contradiction with the definition of $\hat{f}$.

If $\ell$ is infinite, the previous proof still works: one first chooses a very large $l \!\in\! \mathbb{N}$ and approaches $f$ in $C^l$ topology 
by a $C^{\infty}$ diffeomorphism with finitely many fixed points. To do this, one reproduces the argument above but using $C^{\infty}$ 
(and not just $C^l$) diffeomorphisms on maximal intervals $I$ where $f$ is $\eps/2$-close to $\mathrm{id}$ in $C^l$-norm. Since 
the $C^{\infty}$ topology is induced by the family of $C^{l}$-norms, this shows the density of  $G_{1,\infty}$ in $\Diff_P^{\infty}([0,1])$.

Having proved the density of $G_{1,\ell}$, we may assume that the element $f \in \Diff_P^{\ell}([0,1])$ we want to approximate by elements of 
$\mathcal{U}_{\ell}$ has finitely many fixed points. Now, if one of the interior fixed points $p$ is not topologically transversal, we can easily 
destroy it by a $C^{\ell}$-small perturbation just by slightly ``pushing the graph of $f$ around this point in the good direction''.
Formally, this corresponds to 
composing $f$ with a diffeomorphism that is very close to $\mathrm{id}$ in $C^{\ell}$-norm, is supported in a small neighborhood of $p$ 
and moves $p$ to the right (resp. left) if the graph of $f$ is above (below) the diagonal close to $p$. Proceeding this way with all interior 
fixed points of $f$ that are non topologically transversal, we obtain a new diffeomorphism $g$ that is $C^{\ell}$-close to the original one 
but for which all interior fixed points are topologically transversal.

Finally, if $V_{\infty} (g) > 0$ then, by the localization property (\ref{e:localization}), the element $g$ belongs to $\mathcal{U}_{\ell}$. Otherwise, 
take any component $[a,b]$ of $g$. We must have $V_{\infty} (g;[a,b]) = 0$, but as in the proof of Proposition \ref{p:almost-trivial}, any   
nontrivial perturbation $\hat{g}$ of $g$ supported in a single fundamental domain of $g$ inside $(a,b)$ must satisfy $V_{\infty} (\hat{g};[a,b])>0$. 
Thus, this $\hat{g}$ belongs to $\mathcal{U}_{\ell}$, and since $\hat{g}$ can be taken arbitrarily close to $g$ in $C^{\ell}$ topology, 
this concludes the proof.
\end{proof}

\medskip

\begin{rem} 
Although Theorem \ref{t:higher} was stated in regularity $C^{2}$ and higher, a careful reading of the arguments above 
(plus the results from \cite{EN21} and \cite{EN24}) shows that it still holds in regularity $C^{1+bv}$. A crucial 
point is to use the asymptotic variation rather than the Mather invariant, since it is easier to deal with the 
former in this critical regularity.
\end{rem}

\begin{rem}\label{r:mas-arriba}
For $\ell \geq 3$, giving examples of open dense subsets in $\mathrm{Diff}^{\ell}_+ ([0,1])$ in which no element is $C^2$-conjugated to 
a nontrivial power is easier than what was done to prove Theorem \ref{t:higher}. This was somehow discussed in \cite[Section 6]{EN23}, 
where the obstruction is related to the (first) residue at the origin which, for $C^3$ diffeomorphisms, is closely related to the Schwarzian 
derivative (see also \cite{BW04} for results of this type in the real-analytic setting under some extra hypothesis). However, this argument 
breaks down for $\ell = 2$. More importantly, it does not give any further information in relation to distortion elements. 
This is why a global argument (in this case, related to the asymptotic variation) and 
not just a local one is necessary to establish Theorem \ref{t:higher} in full generality.
\end{rem}


\section{A proof of Theorem \ref{t:higher2}}
\label{s:ln}

First recall that, by  the (homogeneity) Property 1, the set of  $C^2$ diffeomorphisms of $[0,1]$ that are $C^2$-conjugated to all  
nontrivial powers of themselves is contained in $\Diff^2_V ([0,1])$. The proof of its density is done in several steps, following a 
strategy similar to that of the proof of Theorem~\ref{t:main}. We fix $f \in \Diff^2_V ([0,1])$.

\vspace{0.25cm}

\noindent{\bf{Step 1.}} Passing to maps with finitely many fixed points.

\medskip

We claim that $f$ is arbitrarily $C^2$-close to elements in $\Diff^2_V ([0,1])$ with finitely many fixed points. 
Indeed, by Property 3, every element in $\Diff_V^2([0,1])$ has only parabolic fixed points. We are hence 
in an analog situation to that of Lemma \ref{l:G1}, except that here we are working in class $C^2$ and not just $C^1$. 
However, an argument in higher regularity was already given at the beginning of the proof of Proposition \ref{p:D2}. 
(Note that this does not use any information on the asymptotic variation.) The only issue is to ensure the vanishing 
of the asymptotic variation for the perturbed map. However, the localization formula allows doing this: it 
suffices, when changing the diffeomorphism on intervals where it is close to the identity, 
to use diffeomorphisms that are the time-1 maps of smooth vector fields.

We will denote by $f_1$ the element in $\Diff^2_V([0,1])$ 
thus obtained that is $C^2$-close to $f$ and has finitely many fixed points.

\medskip

\noindent{\bf{Step 2.}} Changing the germs at the fixed points.

\medskip

We claim that $f_1$ is arbitrarily $C^2$-close to elements in $\Diff^2_p([0,1])$ that,  
about each fixed point, are conjugate by a translation to a germ of the form 
$q_{i}^{\lambda}$ (with $i \in \{1,2,3\}$), where  
\begin{equation}\label{eq:germenes}
q_{1}^{\lambda} (x) := \frac{x}{1 - \lambda x},  
\quad \qquad q_2^{\lambda} (x) := \frac{x}{\sqrt{1-\lambda x^2}} 
\qquad \mbox{ and } \qquad 
q_3^{\lambda} (x) := \frac{x}{\sqrt[3]{1-\lambda x^3}}. 
\end{equation} 
This is similar to Lemma \ref{l:G2}, except that here we need to keep the same second derivative 
at the fixed points along the interpolation. (This is why we need the extra parameter $\lambda$.)  
The precise argument works as follows: Fix a very small $\varepsilon > 0$ so that the $\varepsilon$-neighborhoods 
of the fixed points of $f_1$ are disjoint. For each such a point $a$, denote by $q_{i,a}^{\lambda}$ the conjugate of 
$q_{i}^{\lambda}$ by the translation $T_a$. Now, using bump functions $\rho = \rho_{\varepsilon}$ 
as those of the proof of Lemma \ref{l:G2} but satisfying the extra (and plausible) hypothesis 
$\|D^2 \rho_{\varepsilon}\|_{\infty} \leq 10 / \varepsilon^2$, perform the following procedure:
 
 \vspace{0.2cm}

\noindent - If $D^2 f_1 (a) \neq 0$ then, on $[a-\varepsilon,a+\varepsilon]$, let $f_2$ be equal to $\rho \, q_{1,a}^{\lambda} + (1-\rho) f_1$, 
where $\lambda = D^2 f_1 (a) / 2$. 
 
 \vspace{0.2cm}
 
\noindent -  If $D^2 f_1 (a) = 0$ and $a$ is a topologically transversal fixed point of $f_1$ then, on $[a-\varepsilon,a+\varepsilon]$, 
let $f_2$ be equal to $\rho \, q_{2,a}^1 + (1-\rho) f_1$ if the graph of $f$ is below the diagonal on the left of $a$ (and above on 
the right), and let $f_2 := \rho \, q_{2,a}^{-1} + (1-\rho) f_1$ otherwise.

\vspace{0.2cm}

\noindent -  If $D^2 f_1 (a) = 0$ and $a$ is not a topologically transversal fixed point of $f_1$ then, on $[a-\varepsilon,a+\varepsilon]$, 
let $f_2$ be equal to $\rho \, q_{3,a}^1 + (1-\rho) f_1$ if the graph of $f_1$ is above the diagonal about $a$, and let 
$f_2 := \rho \, q_{3,a}^{-1} + (1-\rho) f_1$ otherwise.

\medskip

We claim that this new map $f_2$ (which depends on the parameter $\varepsilon$) 
is $C^2$-close to $f_1$ provided $\varepsilon$ is properly chosen. To see this, 
denoting $q := q_{i,a}^{\lambda}$, we compute:
$$D^2 f_2 - D^2 f_1 
= 
\rho \, (D^2 q - D^2 f_1) + 2 D\rho \, (Dq - Df_1) + D^2 \rho \, (q - f_1).$$ 
As in the proof of Lemma \ref{l:G2}, the first two terms in the right converge to zero as $\varepsilon$ goes to $0$. To deal with the 
third one note that, by the  the choice of $i$ and $\lambda$ in each case, about each $a$ we have $(q - f_1)(x) = o((x-a)^2)$. 
Therefore, in places where $q - f_1$ is nonzero, one has
$$ |(q- f_1) \, D^2 \rho| \leq \frac{10}{\varepsilon^2} \, |q-f_1| = o(1).$$

Recall that the germs of $f_2$ at the fixed points, being of the form $q_i^{\lambda}$, 
can be conjugated to all their nontrivial powers using the appropriate (germs of) homotheties. 
Note also that $f_2$ depends on $\varepsilon$. Below we will denote $f_{2,\varepsilon}$ this map 
in order to emphasize this dependence.

\medskip

\noindent{\bf{Step 3.}} Forcing the asymptotic variation to vanish again.

\medskip

The asymptotic variation of $f_{2,\eps}$ may have become positive. However, (the continuity) 
Property~4 implies that it must be very small.\footnote{Actually, to show this fact, we do not need to use Property 4 
(which is somewhat hard to prove): it suffices to note that the function  $f \mapsto V_{\infty} (f)$ 
is upper semicontinuous (this easily follows from the definition), hence $0$ is a continuity value of  its image.} 
Morally, this means that the Mather invariant of the diffeomorphism induced by $f_{2,\eps}$ in each of its components is 
``close'' to being trivial, hence it can be destroyed by a small perturbation and, this way, the asymptotic variation 
becomes zero again. Nevertheless, to properly apply this argument in the $C^2$ topology (and not just $C^{1+bv}$), 
we need a slightly more sophisticated argument. 

Since the perturbation we will perform will be supported away from 
the (finitely many) fixed points, with no loss of generality we can 
assume that $f_1$ (and thus $f_{2,\eps}$) has no interior fixed point. 

Let us fix a point $p \in (0,1)$, say $p := \frac12$ for concreteness, and  
for each $\varepsilon > 0$ let us consider $M_{\varepsilon}^{p,p}$, the (representative  of the) Mather invariant of 
$f_{2,\varepsilon}$. By Yoccoz' continuity theorem (cf. \cite[Chapter V]{yoccoz}), as $\eps$ goes to $0$,
this $C^2$ circle diffeomorphism $M_{\varepsilon}^{p,p}$ converges in $C^2$ topology to the 
Mather invariant of $f_1$, which is the identity since $f_1$ has vanishing asymptotic variation (cf. Property 3).

Let us fix $\delta>0$. We claim that, taking $\eps$ small enough (assuming in particular that the intervals $[0,\eps]$, $[f_1^{-3}(p),p]$ 
and $[1-\eps,1]$ are two-by-two disjoint and that $f_{2,\eps}$ is $\delta$-close to $f_1$ in $C^2$ norm), one can construct $f_3$ such that:
\begin{itemize}
\item it coincides with $f_{2,\eps}$ on the $\eps$-neighborhoods of $0$ and $1$;  
\item it coincides with $f_1$ outside these neighborhoods except for two fundamental domains of $f_1$ lying in $[f_1^{-3}(p),p]$; 
\item it has a trivial Mather invariant;
\item it is $\delta$-close to $f_1$ in $C^2$ norm. 
\end{itemize}
To do this, first fix $\eta>0$ such that, if $h_1$ and $h_2$ are any $C^2$ diffeomorphisms 
of $[0,1]$ that are $\eta$-close to the identity in $C^2$ norm, then $f_1\circ h_1\circ h_2$ is $\delta$-close to $f_1$ in $C^2$ norm. 
Using a procedure similar to that of \S \ref{s:cancel}, we will build diffeomorphisms $h_{1,\eps}$ and $h_{2,\eps}$ as above 
in such a way that the diffeomorphism $f_3 := f_{2,\eps}\circ h_{1,\eps}\circ h_{2,\eps}$ satisfies all the requires properties. 

The inverse $(M_{\varepsilon}^{p,p})^{-1}$ of $M_{\varepsilon}^{p,p}$ is a circle diffeomorphism that fixes 
$0$. It can hence be written as a composition $\varphi_{1,\eps}\circ\varphi_{2,\eps}$ of two circle diffeomorphisms supported outside 
$0$ and $\frac12$, respectively, that converge to the identity in $C^2$ topology as $\eps$ goes to $0$. Let $\phi_{1,\eps}$ (resp. 
$\phi_{2,\eps}$) be the diffeomorphism of $\R$ that coincides with a lift of $\varphi_{1,\eps}$ (resp. $\varphi_{2,\eps}$) on $(-1,0)$ (resp. 
$(-\frac52,-\frac32)$) and with the identity elsewhere. Let $(f_t)$ be the $C^1$-flow of which $f_1$ is the time-$1$ map. 
Consider the $C^2$-diffeomorphism $\psi \!: \R \to (0,1)$ given by $\psi(t) := f_t(p)\in (0,1)$ and, for $i \in \{1,2\}$, let 
$h_{i,\eps} := \psi\phi_{i,\eps}\psi^{-1}$. Since $\phi_{1,\eps}$ and $\phi_{2,\eps}$ are supported 
in $[-3,0]$ and $\psi$ is a fixed $C^2$ diffeomorphism sending $[-3,0]$ onto $[f_1^{-3}(p),p]$, the parameter $\eps$ can be chosen 
small enough so that the maps $h_{i,\eps}$ are $\eta$-close to the identity in $C^2$ norm (and supported in $[f^{-3}(p),p]$). 
Now, by Lemma \ref{l:modify}, the diffeomorphism $f_3 := f_{2,\eps}\circ h_{1,\eps}\circ h_{2,\eps}$ has a trivial Mather invariant. 
Note that $f_3$ coincides with $f_1\circ h_{1,\eps}\circ h_{2,\eps}$ on $(\eps,1-\eps)$; hence, by the choice of $\eta$, it is 
$\delta$-close to $f_1$ in $C^2$ norm there. The latter also holds on $[0,\eps]\cup[1-\eps,1]$, since therein $f_3$ coincides 
with $f_{2,\eps}$, which is itself $\delta$-close to $f_1$ in $C^2$ norm. Therefore, $f_3$ has all the required properties.

\medskip

\noindent{\bf{Step 4.}}  Conjugating to nontrivial powers.

\medskip

Finally, it follows as a direct application of Proposition \ref{l:criterion} (with $\ell = 2$) that $f_3$ is $C^2$-conjugated 
to all of its nontrivial powers. This closes the proof of Theorem \ref{t:higher2}.

\vspace{0.2cm}

\begin{rem}\label{r:last}
It is worth noting that the obstruction to being conjugate to nontrivial powers related to residues alluded to in Remark \ref{r:mas-arriba}
persists for diffeomorphisms (with no interior fixed points and) with trivial Mather invariant (hence with vanishing asymptotic variation). 
This implies that, as stated, Theorem \ref{t:higher2} does not extend to $\Diff^{\ell}_V ([0,1])$ for $\ell \geq 3$. However, it should not 
be difficult to establish an analog version of it adding the hypothesis of vanishing of the residues at the fixed points (whenever 
they are defined). We leave this to the interested reader.
\end{rem}

\begin{rem} 
In contrast to Theorem \ref{t:higher2}, it is worth pointing out that the subset of $\Diff_V^2 ([0,1])$ formed by the diffeomorphisms 
that are not $C^2$-conjugated to any nontrivial power of themselves is also dense in $\Diff_V^2 ([0,1])$. The proof of this is exactly 
the same except that one needs to replace (\ref{eq:germenes}) by germs that are not $C^2$-conjugated to their nontrivial powers. 
In concrete terms, one replaces $q_1^{\lambda}$ by 
$$\hat{q}^{\lambda}_1 (x) := x + \lambda x^2$$
which, according to \cite{EN24}, is not $C^2$-conjugated to any of its nontrivial powers (because of the residue). Concerning 
$q_2^{\lambda}$ and $q_3^{\lambda}$, one may replace them by any infinitely-flat germ with trivial $C^2$-centralizer 
and/or its inverse on each side of the fixed point (according to the type of transversality). As a concrete example, one 
may use the germ constructed by Sergeraert in \cite{Se77}, which was shown to have trivial $C^2$-centralizer in \cite{Ey09}. 
\end{rem}


\section{Some final remarks}
\label{s:remarks-fin}

There are other sources of distortion elements that are not related to conjugacy to nontrivial powers, particularly in low regularity. 
For instance, recall that every finitely-generated torsion-free nilpotent group realizes as a group of $C^1$ diffeomorphisms of the interval 
\cite{FF03, CJN14,Pa16}. Now, in a torsion-free nilpotent group that is not almost Abelian, elements in the center turn out to be distorted, yet they 
are not conjugated to nontrivial powers of themselves \cite[Appendix]{Gr81}.

As a concrete example, consider the Heisenberg group 
$$H = \big\langle f,g,h : [f,g]=h, [f,h]=[g,h]= \mathrm{id} \big\rangle.$$ 
An easy computation 
shows that $[f^n,g^n] = h^{n^2}$, from where it follows that $h$ is distorted. The realization of $H$ as a group of diffeomorphisms 
of the interval described in \cite{FF03,CJN14,Pa16} actually embeds it into $\Diff^1_P ([0,1])$. This way, $h$ becomes a diffeomorphism of the 
interval that is $C^1$-distorted. It is worth pointing out that $h$ has infinitely many components (actually, this holds for any faithful action of 
$H$ by {\em homeomorphisms} of the interval \cite[Section 2.2.5]{Na11}). The argument of the first proof of Theorem~\ref{t:generic} then 
shows that $h$ is not $C^1$-conjugated to a nontrivial power of itself.

\medskip

\noindent{\bf Question 4.} Can $H$ be realized as a subgroup of $\Diff^1_P([0,1])$ so that the element $h$ becomes a $C^2$ diffeomorphism~? 
If so, can $V_{\infty}(h)$ take arbitrary (nonegative) values~?

\medskip

Recall that, by the Plante-Thurston's theorem \cite{PT76}, the whole group $H$ cannot 
be realized as a group of $C^{1+bv}$ diffeomorphisms of the interval 
(see \cite[Section 4.2]{Na11} for a full discussion on this).


\vspace{0.4cm}

\noindent{\bf Acknowledgments.} Both authors are very grateful to Christian Bonatti for having 
shared some of his views (including an ``old'' manuscript) related to a crucial section of this paper, 
as well as to all the organizers of the special trimester ``Group actions and rigidity: around the Zimmer 
program'' for providing us the perfect conditions to write up this article. 

This work was supported by the ECOS project 23003 ``Small spaces under actions''. The authors 
also acknowledge support of the Institut Henri Poincar\'e (UAR 839 CNRS-Sorbonne Universit\'e), 
and LabEx CARMIN (ANR-10-LABX-59-01).


\begin{footnotesize}

\vspace{0.25cm}

\noindent {\bf H\'el\`ene Eynard-Bontemps} \hfill{\bf Andr\'es Navas}

\noindent Universit\'e Grenoble Alpes, CNRS\hfill{ Dpto. de Matem\'aticas y C.C.}

\noindent  Institut Fourier  \hfill{ Universidad de Santiago de  Chile}

\noindent 38000 Grenoble \hfill{Alameda Bernardo O'Higgins 3363}

\noindent France \hfill{Estaci\'on Central, Santiago, Chile} 

\noindent helene.eynard-bontemps@univ-grenoble-alpes.fr \hfill{andres.navas@usach.cl}

\end{footnotesize}

\end{document}